\documentclass[journal]{IEEEtran}
\usepackage{epsf, psfrag, amssymb, amsfonts, amsmath,enumerate,cite}
\usepackage{graphicx, subfigure, color,bbm, bm}
\usepackage{boxedminipage}
\usepackage{multirow}
\usepackage{booktabs}
\usepackage{algorithm, algorithmic}

\makeatletter
\def\hlinewd#1{%
  \noalign{\ifnum0=`}\fi\hrule \@height #1 \futurelet
   \reserved@a\@xhline}
\makeatother

\floatstyle{plain}
\newfloat{myalgo}{tbhp}{mya}
\newenvironment{Algorithm}[2][tbh]%
{\begin{myalgo}[#1]
\centering
\begin{minipage}{#2}
\begin{algorithm}[H]}%
{\end{algorithm}
\end{minipage}
\end{myalgo}}
\newtheorem{remark}{Remark}
\usepackage[scriptsize,bf]{caption}
\usepackage{etoolbox}
\AtBeginEnvironment{align}{\setcounter{subeqn}{0}}
\newcounter{subeqn} %

\newtheorem{theorem}{Theorem}
\newtheorem{example}{Example}
\newtheorem{corollary}{Corollary}
\newtheorem{lemma}{Lemma}
\newtheorem{definition}{Definition}

\def\psfancypar#1#2{\begingroup\def\par{\endgraf\endgroup\lineskiplimit=0pt}
               \setbox2=\hbox{\large\sc #2}
               \newdimen\tmpht \tmpht \ht2 \advance\tmpht by \baselineskip
               \font\hhuge=Times-Bold at \tmpht
               \setbox1=\hbox{{\hhuge #1}}
               \count7=\tmpht \count8=\ht1
               \divide\count8 by 1000 \divide\count7 by \count8
               \tmpht=.001\tmpht\multiply\tmpht by \count7
               \font\hhuge=Times-Bold at \tmpht
               \setbox1=\hbox{{\hhuge #1}}
               \noindent
                \hangindent1.05\wd1
               \hangafter=-2 {\hskip-\hangindent
               \lower1\ht1\hbox{\raise1.0\ht2\copy1}%
                \kern-0\wd1}\copy2\lineskiplimit=-1000pt}

\newcommand{\beq}{\begin{equation}}
\newcommand{\eeq}{\end{equation}}
\newcommand{\bqa}{\begin{eqnarray}}
\newcommand{\eqa}{\end{eqnarray}}
\newcommand{\bqn}{\begin{eqnarray*}}
\newcommand{\eqn}{\end{eqnarray*}}
\newcommand{\nn}{\nonumber}

\newcommand{\be}{\begin{enumerate}}
\newcommand{\ee}{\end{enumerate}}
\newcommand{\bi}{\begin{itemize}}
\newcommand{\ei}{\end{itemize}}
\newcommand{\bd}{\begin{description}}
\newcommand{\ed}{\end{description}}
\newcommand{\ba}{\begin{array}}
\newcommand{\ea}{\end{array}}
\newcommand{\bde}{\begin{definition}}
\newcommand{\ede}{\end{definition}}
\newcommand{\bex}{\begin{example}}
\newcommand{\eex}{\end{example}}

\def\boxit#1{\vbox{\hrule\hbox{\vrule\kern3pt
        \vbox{\kern3pt#1\kern3pt}\kern3pt\vrule}\hrule}}

\def\reals{ { {\rm  I \kern-0.15em R }  } }
\def\complex{ {\,{{\rm C} \kern-0.50em \raise0.20ex {  |}}\, }}

\def\0bf{{\bf 0}}
\def\1bf{{\bf 1}}
\def\2bf{{\bf 2}}
\def\3bf{{\bf 3}}
\def\4bf{{\bf 4}}
\def\5bf{{\bf 5}}
\def\6bf{{\bf 6}}
\def\7bf{{\bf 7}}
\def\8bf{{\bf 8}}
\def\9bf{{\bf 9}}

\def\Rbf{{\bf R}}

\def\Rxx{\Rbf_{\ssstyle X\kern-.1em X}}

\let\ssstyle=\scriptscriptstyle

\def\Kout{\setbox1=\hbox{\Huge\bf K}\hbox to
1.05\wd1{\hspace{.05\wd1}
}}
\def\Sout{\setbox1=\hbox{\Huge\bf S}\hbox to 1.05\wd1{\hspace{.05\wd1}
}}

\DeclareMathOperator*{\argmin}{arg\,min}

\begin{document}

\title{Distributed Average Consensus with Bounded Quantizer and Unbounded Input}
\author{Shengyu Zhu and Biao Chen
\thanks{This work was presented in part at the IEEE International Workshop on Signal Processing  Advances in Wireless Communications (SPAWC), Edinburgh, UK, July 2016 \cite{Zhu2016SPAWC}.}
\thanks{The authors are with the Department of Electrical Engineering and Computer Science, Syracuse University, Syracuse, NY 13244, USA (e-mail: szhu05@syr.edu; bichen@syr.edu).}}

\maketitle

\begin{abstract}
This paper considers distributed average consensus using finite-bit bounded quantizer with possibly unbounded data. Under the framework of the alternating direction method of multipliers (ADMM), we develop distributed averaging algorithms where each node iteratively updates using only the local information and finitely quantized outputs from its neighbors. It is shown that all the agent variables either converge to the same quantization level or cycle around the data average after finite iterations. An error bound for the consensus value is established, which turns out to be the same as that of using the unbounded rounding quantizer provided that an algorithm parameter (i.e., ADMM step size) is small enough. We also analyze the effect of the algorithm parameter and propose an adaptive parameter selection strategy that only requires knowledge of the number of agents in order to accelerate the algorithm with certain consensus accuracy guarantee. Finally, simulations are performed to illustrate the effectiveness of the proposed algorithms.
\end{abstract}

\section{Introduction}
\IEEEPARstart{T}{here} have been extensive efforts devoted to applications where autonomous agents collaborate to accomplish a global objective. Of particular interest is distributed average consensus which computes the average of agents' measurements using only local computation and communication. The problem originates from distributed computation and decision-making \cite{Tsitsiklis1984,Tsitsiklis1986}, and arises in various applications including load balancing for parallel computers \cite{Xu1996}, multi-agent coordination \cite{Ren2005}, gossip algorithms \cite{Boyd2006}, and distributed inference \cite{Kar2008}. We refer the reader to \cite{Cao2013} for an overview.

Distributed averaging algorithms are iterative algorithms that aim to achieve the average consensus for all agents. They are more energy efficient and are robust to link failures compared with fusion center based processing. Studied in \cite{Elsner1990, Xiao2004, Nedic2009} is a classical method where each agent updates itself with a weighted average of its own value and the values received from its neighbors. Another approach is a gossip based algorithm, initially introduced in  \cite{Tsitsiklis1984} for consensus problems and further investigated in \cite{Aysal2009,Boyd2006}, among others. In \cite{Schizas2008,Zhu2009,Erseghe2011}, the authors proposed a distributed averaging algorithm based on the alternating direction method of multipliers (ADMM), an iterative algorithm for solving convex problems and has received much attention recently (see \cite{BoydADMM}). The idea is to formulate the data average as the solution to a least-squares problem and then decentralize the ADMM update. While the above approaches converge to the desired average for each agent under fairly general conditions, it is important to note that agents are assumed to communicate real values of infinite precision with their neighbors. In practical applications such as those encountered in wireless sensor networks, data communications are usually subject to limited capacity and resource constraints. As such, it is commonly assumed that agents can reliably transmit only quantized data, resulting in {\emph{quantized consensus}} \cite{Kashyap2007}.

{\em Literature review on quantized consensus:} A well studied approach to quantized consensus is to use dithered quantizers \cite{Schuchman1964}. Utilizing the first- and second-order moments of dithered quantizer output, many researchers have developed algorithms that yield consensuses at a random variable whose mean is the desired average, along with analyses on the mean squared error and convergence rate \cite{Aysal2008,Zhu2015,Kar2010,Carli2010,Zhu2009,Erseghe2011}. Deterministic algorithms for quantized consensus have also been studied. In \cite{Nedic2009,Carli2010}, the proposed algorithms ensure the convergence to a consensus in a neighborhood of the average within finite iterations and consensus accuracy is characterized by an error bound on the difference between the consensus result and the desired average. Recent work of \cite{Zhu2016TSP,Zhu2016TSPerr} reports a similar result for an ADMM based distributed averaging algorithm using rounding quantizer, which is referred to as deterministically quantized consensus ADMM (DQ-CADMM). An advantage of DQ-CADMM is that the error bound does not depend on the size of the network nor the agents' data. The consensus algorithm in \cite{Chamie2014} is shown to either reach a quantized consensus close to the average in finite time or lead all agent variables to cycle in a small neighborhood around the average; in the latter case, consensus is not guaranteed. In \cite{Li2011}, quantized consensus is formulated as a feedback control design problem for coding/decoding schemes. With an appropriate scaling function and carefully chosen control gain based on some spectral properties of the Laplacian matrix of the underlying fixed undirected graph, the proposed protocol with rounding quantizer can achieve the exact average consensus asymptotically. 

We note that the quantizers in most existing works are still of infinite bits since the quantizer output has an unbounded range and infinite quantization levels. In order to achieve finite-bit communications at each iteration, these algorithms require the knowledge of a bound on agents' data such that truncation can be employed. It is worth mentioning that there also exist algorithms that can work with predefined finite-bit quantizers; instead of directly truncating the communicated data, the bound on agents' data is used to pick appropriate updating weights to handle the finite-bit communication constraint. With the same bounded condition, the consensus algorithms in \cite{Kar2010} and \cite{Li2011} can use finite-level quantizers to asymptotically achieve the average for each agent.  
{\it Summary of contributions:} We study distributed average consensus using finite-bit {\em bounded} quantized communications with possibly {\em unbounded} inputs, i.e., without knowing a bound on agents' data. This is primarily useful to statistical inference in sensor networks where the observations are intrinsically unbounded, e.g., when observations are corrupted by additive white Gaussian noise (AWGN). Naively truncating agents' data to meet the bounded condition may lead to loss of inference performance. For example, considered in \cite{Zhu2016ISIT, Zhu2016Detection} is distributed detection using consensus based approach where agents' data for averaging are the log-likelihood ratios of local observations. Without properly handling the unbounded data, existing approaches fail to establish the optimal asymptotic detection performance using finite-bit information exchange per iteration. Another example is from parameter estimation. With AWGN, it is well known that the global average of local observations is a minimum variance unbiased estimator for the parameter. Thus,  the data for averaging in consensus based approaches are the local observations at each node and are again unbounded. If truncation is used to make local data bounded, consensus based approach may result in biased estimator.

In this paper, we consider uniform deterministic quantization that is finite-bit and bounded. Our quantizer is constructed by projecting the input to a compact convex set and then applying the uniform rounding quantizer to the projected value. We propose an ADMM based quantized consensus algorithm by replacing the rounding quantizer in DQ-CADMM \cite{Zhu2016TSP,Zhu2016TSPerr} with this bounded quantizer. Since the projection operation may lead to unbounded quantization error with a bounded quantizer, the convergence analysis becomes more complicated. To proceed, we show that the proposed algorithm is equivalent to applying the usual rounding quantizer to the ADMM update of a constrained least-squares problem and that the agent variables are eventually bounded. We establish that the algorithm, within finite iterations, either converges to a common quantization point or cycles with a finite period around the average; in the latter case, every agent has the same sample average of quantized variable values over one period. Furthermore, using convexity properties of the constrained least-squares problem, we derive an upper bound on the consensus error which does not depend on the size of the network nor agents' data. We also discuss the effect of the ADMM step size on the performance of the proposed algorithm and come up with a heuristic parameter selection scheme that requires the knowledge of the number of nodes and the number of edges (indeed, it suffices to know only the number of agents or an upper bound on it).  


Preliminary results were reported in \cite{Zhu2016SPAWC}. Besides providing thorough simulations on the proposed algorithm, the current paper characterizes the consensus error in cyclic cases and investigates the effect of the algorithm parameter. An adaptive strategy is then proposed for parameter selection, which is shown to significantly accelerate the algorithm with consensus accuracy guarantee.

{\em Organization of the paper:} The rest of this paper is organized as follows. Section \ref{sec:Pre} reviews consensus ADMM for distributed optimization and its application to distributed average consensus. In Section \ref{sec:BQC}, we define the bounded quantizer and propose a distributed averaging algorithm that uses this quantizer. We establish its convergence result and further come up with an improved algorithm that results in the same error bound as DQ-CADMM. Section \ref{sec:rho} discusses the effect of the algorithm parameter. Simulations are provided in Section \ref{sec:simulation} and Section \ref{sec:conclusion} concludes the paper.

     

{\it Notations:} 
We use $\bm 0$ to denote the all-zero column vector with a suitably defined dimension. Notation $\bm 1_n$ denotes the $n$-dimensional all-one column vector; $\bm 0_n$ and $\bm I_n$ are the $n\times n$ all-zero and identity matrices, respectively. For $\tilde{x}\in\mathbb{R}$, $\lceil \tilde{x}\rceil$ denotes the smallest integer that is greater than or equal to $x$. Notation $\|\bm x\|_2$ denotes the Euclidean norm of a vector $\bm x$. Given a positive semidefinite matrix $\bm G$ with proper dimensions, the $\bm G$-norm of $\bm x$ is $\|\bm x\|_{\bm G}=\sqrt{\bm x^T\bm G\bm x}$. For a symmetric matrix $\bm L\in\mathbb{R}^{n\times n}$, denote its eigenvalues in the ascending order as $\lambda_1(\bm L)\leq  \lambda_2(\bm L)\leq\cdots\leq\lambda_n(\bm L)$. For any matrix $\bm M$, $\mathcal{C}(\bm M)$ represents its column space.


\section{Consensus ADMM for Distributed Averaging}
\label{sec:Pre}
This section briefly reviews consensus ADMM (CADMM) and its application to distributed average consensus. We start with the network model used in this paper.
\subsection{Network Model}
Consider a network of $n$ agents bi-directionally connected by $m$ edges (hence $2m$ arcs), and with a fixed topology. We describe this network as a symmetric directed graph $\mathcal{G}_d=\{\mathcal{V},\mathcal{A}\}$ or an undirected graph $\mathcal{G}_u=\{\mathcal{V},\mathcal{E}\}$, where $\mathcal{V}=\{1,2,\dots,n\}$ is the set of vertices, $\mathcal{A}=\{(i,j):~\text{the arc from agent}~i~\text{to agent}~j\}$ is the set of arcs, and $\mathcal{E}$ is the set of edges with $|\mathcal{E}|=m$. Define the oriented incidence matrix $\bm{M}_-\in\mathbb{R}^{n\times 2m}$ with respect to $\mathcal{G}_d$ as follows: $[\bm M_-]_{i,l}=1$ if the $l$th arc leaves agent $i$, $[\bm M_-]_{i,l}=-1$ if the $l$th enters agent $i$, and $[\bm M_-]_{i,l}=0$ otherwise. The unoriented incidence matrix $\bm{M}_+\in\mathbb{R}^{n\times2m}$ is defined by setting $[\bm M_+]_{i,l}=\left|[\bm M_-]_{i,l}\right|$. Denote $\mathcal{N}_i=\{j:(i,j)\in\mathcal{A}\}$ as the set of neighbors of agent $i$. Further define $\bm L_-=\frac{1}{2}\bm M_-\bm M_-^T$ and $\bm L_+=\frac{1}{2}\bm M_+\bm M_+^T$  which are respectively the signed and signless Laplacian matrices with respect to $\mathcal{G}_u$. Then $\bm W=\frac{1}{2}(\bm L_-+\bm L_+)=\text{diag}\{|\mathcal{N}_1|,|\mathcal{N}_2|,\cdots,|\mathcal{N}_n|\}$ is the degree matrix related to $\mathcal{G}_u$, i.e., a diagonal matrix with $(i,i)$th entry being $|\mathcal{N}_i|$ and other entries being $0$. 

The following lemma is stated to help establish  our main results in Section~\ref{sec:BQC}.
\begin{lemma}[\hspace{-0.1pt}\cite{ChungSpectral,FiedlerAlgebraic}]
\label{lem:graphLemma}
Given a connected network, we have
\begin{enumerate}[~~a)]
\item $\bm L_-$ is positive semidefinite and $0=\lambda_1(\bm L_-)<\lambda_2(\bm L_-)\leq\lambda_3(\bm L_-)\leq\cdots\leq\lambda_n(\bm L_-)$. $\bm L_- \bm b=\bm 0$ if and only if $\bm b\in\mathcal{C}(\bm1_n)$.
\item $\bm L_+$ is positive semidefinite and $\lambda_n(\bm L_+)>0$.
\item $\mathcal{C}(\bm M_-)=\mathcal{C}(\bm L_-)$. For every $\bm\alpha\in\mathcal{C}(\bm L_-)$, there exists a unique $\bm\beta\in\mathcal{C}({\bm M_-^T})$ such that $\bm\alpha=\bm M_-\bm\beta$.
\end{enumerate}
\end{lemma}
\subsection{CADMM for Distributed Average Consensus} 
\label{sc:CADMM_DAC}
Let $r_i\in\mathbb{R}$ be the local data at node $i$, $i=1,2,\ldots,n$. Distributed average consensus computes the global average
$$\bar{r}=\frac{1}{n}\sum_{i=1}^n r_i,$$ 
using only local information exchange among neighboring nodes. We assume that nodes do not know the global topology of the network, which is typical in large scale networks or when nodes are autonomous. As a result, the communicated information among linked nodes at each iteration is usually a real value or its quantized version in many distributed averaging algorithms.

The original CADMM is obtained by decentralizing the ADMM update of minimizing a sum of convex functions. Let $f_i(\tilde{\bm x}):\mathbb{R}^d\to\mathbb{R}$, where $d$ is a positive integer, denote a convex local objective function that is only known to node $i$. It has been shown in \cite{Schizas2008,Shi2014} that the CADMM update for minimizing 
$\sum_{i=1}^nf_i(\tilde{\bm x})$ is
\begin{equation}
\begin{aligned}
\label{eqn:distributedversion}
\bm x_i^{k+1}=&~\argmin_{\tilde{\bm x}} f_i(\tilde{\bm x})+\rho|\mathcal{N}_i|\tilde{\bm x}^T\tilde{\bm x}+\tilde{\bm x}^T\Bigg(\rho|\mathcal{N}_i|\bm x_i^k\\&~+\rho\sum_{j\in\mathcal{N}_i}\bm x_j^k-\bm\alpha_i^k\Bigg),\\
\bm \alpha_i^{k+1}=&~\bm\alpha_i^k+\rho\Bigg(|\mathcal{N}_i|\bm x_i^{k+1}-\sum_{j\in\mathcal{N}_i}\bm x_j^{k+1}\Bigg),
\end{aligned}
\end{equation}
where $\rho>0$ is the ADMM parameter and $\bm\alpha_i^k\in\mathbb{R}^d$ is the local Lagrangian multiplier of node~$i$. The above update is fully decentralized as the update of ${\bm x}_i^{k+1}$ and $\bm\alpha_i^{k+1}$ only relies on local and neighboring information. We refer to (\ref{eqn:distributedversion}) as the original CADMM. 


The original CADMM converges under fairly mild conditions, which follows directly from global convergence of the ADMM \cite{BoydADMM,He2012,Deng2016}. To apply CADMM to distributed average consensus, it is noted the global average is the unique solution to a least-squares problem, i.e., $$\bar{r}=\argmin_{\tilde{x}}\frac{1}{2}\sum_{i=1}^n(\tilde{x}-r_i)^2.$$
The corresponding CADMM update for distributed averaging is then given by
\begin{equation}
\label{eqn:CADMMDA}
\begin{aligned}
x_i^{k+1}&=\frac{1}{1+2\rho|\mathcal{N}_i|}\Bigg(\rho|\mathcal{N}_i|x_i^k+\rho\sum_{j\in\mathcal{N}_i}x_j^k-\alpha_i^k+r_i\Bigg),\\
\alpha_i^{k+1}&=\alpha_i^k+\rho\Bigg(|\mathcal{N}_i|x_i^k-\sum_{j\in\mathcal{N}_i}x_j^k\Bigg).
\end{aligned}
\end{equation}

For ease of presentation, we further write (\ref{eqn:CADMMDA}) in a compact form. Define $\bm x^k$ and $\bm\alpha^k$ as the vectors concatenating all $x_i^k$ and $\alpha_i^k$, respectively. Then (\ref{eqn:CADMMDA}) is equivalent to
\begin{align}
\bm x^{k+1}&=(\bm I_n+2\rho\bm W)^{-1}\left(\rho\bm L_+\bm x^k-\bm\alpha^k+\bm r\right),\label{eqn:QCADMMmatrixa}\\
\bm\alpha^{k+1}&=\bm\alpha^k+\rho\bm L_-\bm x^k.\label{eqn:QCADMMmatrixb}
\end{align}
Letting $\bm s^k=[\bm x^k;\bm\alpha^k;\bm r]$ and $\bm D_0=(\bm I_n+2\rho\bm W)^{-1}$, we have 
\begin{equation}
\label{eqn:iterateform}
\begin{aligned}
\bm s^{k+1}=\bm D\bm s^k=\bm D^{k+1}\bm s^0,
\end{aligned}
\end{equation}
where $\bm D^{k+1}$ denotes the $(k+1)$th power of the square matrix $\bm D$ defined by
\begin{equation}
\begin{aligned}
\bm D \triangleq \begin{bmatrix} \rho \bm D_0\bm L_+ & -\bm D_0& \bm D_0 \\ \rho^2 \bm L_-\bm D_0\bm L_+ & \bm I_n-\rho \bm L_-\bm D_0& \rho \bm L_-\bm D_0 \\\bm 0_n & \bm 0_n & \bm I_n \end{bmatrix}.\nn
\end{aligned}
\end{equation}
The convergence result of CADMM implies that
\begin{align}
\label{eqn:optsk}
\lim_{k\to\infty} \bm s^k  =\bm s^*=\begin{bmatrix}\bm x^*\\ \bm\alpha^*\\\bm r\end{bmatrix}=\begin{bmatrix} \bm1_n\bar{r}\\ \bm r-\bm1_n\bar{r}\\\bm r\end{bmatrix},
\end{align}
provided that $\bm\alpha^0$ is initialized in $\mathcal{C}(\bm L_-)$ (e.g., all $0$). Moreover, $\bm s^k$ converges linearly to $\bm s^*$. Let $\bm z^k=\frac{1}{2}\bm M_+\bm x^k$ and $\bm\beta^k$ be the unique vector in $\mathcal{C}(\bm M_-^T)$~such that $\bm\alpha^k = \bm M_-\bm \beta^k$ (cf. Lemma~\ref{lem:graphLemma}). Define also 
\begin{align}
\bm u^k = \begin{bmatrix} \bm z^k\\\bm \beta^k \end{bmatrix}~\text{and}~\bm G = \begin{bmatrix} \rho \bm I_{2m} & \bm 0_{2m} \\ \bm 0_{2m}& \frac{1}{\rho}\bm I_{2m}  \end{bmatrix}.\nn
\end{align}
The linear convergence of $\bm s^k\to\bm s^*$ is formally stated below.
\begin{theorem}[Linear convergence of CADMM for distributed average consensus{\cite{Shi2014}}]
\label{thm:linearconvergence} 
Let $\bm x^0\in\mathbb{R}^n$ and $\bm\alpha^0\in\mathcal{C}(\bm L_-)$. Then $\bm u^k$ converges Q-linearly to $\bm u^*=[\bm z^*;\bm\beta^*]$, with $\bm z^*=\bm1_n\bar{r}$ and $\bm\beta^*$ being the unique vector in $\mathcal{C}(\bm M_-)$ such that $\bm M_-^T\bm\beta^*=\bm r-\bm1_n\bar{r}$, with respect to the $\bm G$-norm:
\begin{equation}
\begin{aligned}
\|\bm u^{k+1}-\bm u^*\|^2_{\bm G}\leq\frac{1}{1+\delta}\|\bm u^k-\bm u^*\|^2_{\bm G},\nn
\end{aligned}
\end{equation}
where $\mu>1$ can be any and 
\begin{align}
\delta &= \min \left\{\frac{(\mu-1)\lambda_2(\bm L_-)}{\mu \lambda_n(\bm L_+)},\frac{2\rho \lambda_2(\bm L_-)}{\rho^2\lambda_n(\bm L_+)\lambda_2(\bm L_-)+\mu}\right\}.\nn
\end{align}
Furthermore, $\bm s^k$ is R-linearly convergent to $\bm s^*$ as 
\begin{align}
\|\bm s^{k+1}-\bm s^*\|_2\leq\left(1+\sqrt{\frac{2\rho\lambda_n(\bm L_-)}{1+\delta}}\right)\bm\|\bm u^k-\bm u^*\|_{\bm G}.\nn
\end{align}
\end{theorem}

\section{Distributed Average Consensus with Finite-Bit Bounded Quantization}
\label{sec:BQC}
This section studies distributed average consensus subject to bounded quantization constraint without a known  bound on agents' data. Our finite-bit quantizer is defined as follows. 
\subsection{Finite-bit Bounded Quantizer}
We assume that each agent can store and compute real values with infinite precision but can only send quantized data through the channel which are received by its neighbors without any error. With a predefined quantization resolution $\Delta>0$, the quantization lattice in $\mathbb{R}$ is $$ \Lambda = \{t\Delta: t\in\mathbb{Z}\}.$$ Consider the usual rounding quantizer $\mathcal{Q}: \mathbb{R}\to\Lambda$ that maps a real value to the nearest point in $\Lambda$: 
\begin{equation}
\begin{aligned}
\mathcal{Q}(\tilde{x}) = t \Delta,~\text{if}~\left(t-\frac{1}{2}\right)\Delta< \tilde{x}\leq \left(t+\frac{1}{2}\right)\Delta.\nn
\end{aligned}
\end{equation}
Let $e(\tilde{x})=\mathcal{Q}(\tilde{x})-\tilde{x}$ denote the quantization error. It is clear that $|e(\tilde{x})|\leq\frac{1}{2}\Delta$
for any $\tilde{x}\in\mathbb{R}$. Noticing that the rounding quantizer has infinite quantization levels due to its unbounded range, we further define a bounded quantizer as follows.

Let $\mathcal{X}$ be a nonempty compact convex set in $\mathbb{R}$. We assume without loss of generality that $\mathcal{X}=[-L,L]$ for some $0<L<\infty$, since we can always translate the set. Further assume that $L$ is a multiple of $\Delta$ and let $\Lambda_\mathcal{X}=[-L,L]\cap\Lambda$. Define $\mathcal{T_X}(\cdot)$ as the projection operator that maps a real value to the nearest point in $\mathcal{X}$. A bounded quantizer ${\mathcal{Q}_b}:\mathbb{R}\to\Lambda_\mathcal{X}$ is defined by first projecting its argument onto $\mathcal{X}$ and then applying the usual rounding quantizer $\mathcal{Q}(\cdot)$ to the projected value, i.e.,
\begin{align}
\mathcal{Q}_b(\cdot)=\mathcal{Q}\circ\mathcal{T}_\mathcal{X}(\cdot).\nn
\end{align} 

It is straightforward to see that $|\mathcal{Q}_b(\tilde{x})|\leq L$ for any $\tilde{x}\in\mathbb{R}$. Therefore, the bounded quantizer $\mathcal{Q}_b(\cdot)$ has $2{L}/{\Delta}+1$ quantization levels and its output can be represented by $\lceil\log_2(2{L}/{\Delta}+1)\rceil$ bits. Define $e_b(\tilde{x})=\mathcal{Q}_b(\tilde{x})-\tilde{x}$ to be the quantization error of $\mathcal{Q}_b(\cdot)$, then $|e_b(\tilde{x})|$ can be unbounded for an unbounded $\tilde{x}$. As a note, the above operators operate on each entry of the argument when they have a vector input.

\subsection{CADMM with Finite-bit Bounded Quantization}

We now modify the original CADMM update in (\ref{eqn:CADMMDA}) using the bounded quantizer $\mathcal{Q}_b(\cdot)$, 
as presented in Algorithm \ref{tab:QCADMM}, which we refer to as bounded quantizer based CADMM (BQ-CADMM).

\begin{Algorithm}{\linewidth}
	\caption{BQ-CADMM for distributed average consensus}
	\begin{algorithmic}[1]
	\label{tab:QCADMM}
	\REQUIRE Initialize~$x_i^0=0$ and $\alpha_{i}^0=0$ for each agent $i,i=1,2,\ldots,n$. Set $\rho>0$ and $k=0$.
	\REPEAT
			\STATE every agent $i$ {\bf do}
			\begin{align}
			x_i^{k+1}=&~\frac{1}{1+2\rho|\mathcal{N}_i|}\Bigg(\rho|\mathcal{N}_i|\mathcal{Q}_b(x_i^k)+\rho\sum_{j\in\mathcal{N}_i} \mathcal{Q}_b(x_j^k)\nn\\&~-\alpha_i^k+r_i\Bigg),\nn\\
\alpha_i^{k+1}=&~\alpha_i^k+\rho\Bigg(|\mathcal{N}_i|\mathcal{Q}_b(x_{i}^{k+1})-\sum_{j\in\mathcal{N}_i}\mathcal{Q}_b(x_{j}^{k+1})\Bigg).\nn
			\end{align}
			\STATE {\bf set} $k=k+1$.	
	\UNTIL{a predefined stopping criterion (e.g., a maximum iteration number) is satisfied.}
	\end{algorithmic}
\end{Algorithm}

It is clear that BQ-CADMM is obtained from the CADMM update (\ref{eqn:CADMMDA}) by applying the bounded quantizer $\mathcal{Q}_b(\cdot)$ to local variables $x_i^k$. Here we  provide an alternative interpretation of this algorithm.  Since $\mathcal{Q}_b(\cdot)=\mathcal{Q}\circ\mathcal{T_X}(\cdot)$, BQ-CADMM can be obtained by applying the rounding quantizer $\mathcal{Q}(\cdot)$ to $\mathcal{T_X}(x_i^k)$ in the following update:
\begin{align}
x_i^{k+1}=&~\frac{1}{1\hspace{-1pt}+\hspace{-1pt}2\rho|\mathcal{N}_i|}\Bigg(\rho|\mathcal{N}_i|\mathcal{T_X}(x_i^k)+\rho\sum_{j\in\mathcal{N}_i}\mathcal{T_X}(x_j^k)\nn\\&~-\alpha_i^k+r_i\Bigg),\label{eqn:ProjCADMMa}\\
\alpha_i^{k+1}=&~\alpha_i^k+\rho\Bigg(|\mathcal{N}_i|\mathcal{T_X}(x_i^{k+1})-\sum_{j\in\mathcal{N}_i}\mathcal{T_X}(x_j^{k+1})\Bigg).\label{eqn:ProjCADMMb}
\end{align}
That is, at the $k$th iteration, node $i$ has its local variable value $x_i^k$ and sends its projection $\mathcal{T_X}(x_i^k)$ to the rounding quantizer $\mathcal{Q}(\cdot)$. At the $(k+1)$th iteration, $x_i^{k+1}$ and $\alpha_i^{k+1}$ are updated using  $\mathcal{Q}(\mathcal{T_X}(x_i^k))$ and $\mathcal{Q}(\mathcal{T_X}(x_j^k))$. Since $\mathcal{T_X}(x_i^{k})$ is the projection of $x_i^{k}$ onto $\mathcal{X}$, we can write  (\ref{eqn:ProjCADMMa}) equivalently as
\begin{align}
&\hspace{2pt}\mathcal{T_X}(x_i^{k+1})\nn\\
=&\hspace{2pt}\argmin_{\tilde{x}\in\mathcal{X}}\frac{1}{2}(\tilde{x}-r_i)^2+\rho|\mathcal{N}_i|\tilde{x}^2+\tilde{x}\Bigg(\rho|\mathcal{N}_i|\mathcal{T_X}(x_i^k)\nn\\&+\rho\sum_{j\in\mathcal{N}_i}\mathcal{T_X}(x_j^k)-\alpha_i^k\Bigg)\nn\\
\label{eqn:ConstrainedDOpt}
=&\hspace{2pt}\argmin_{\tilde{x}}\frac{1}{2}(\tilde{x}-r_i)^2+I_\mathcal{X}(\tilde{x})+\rho|\mathcal{N}_i|\tilde{x}^2-\tilde{x}\Bigg(\rho|\mathcal{N}_i|\mathcal{T_X}(x_i^k)\nn\\
&\hspace{2pt}+\rho\sum_{j\in\mathcal{N}_i}\mathcal{T_X}(x_j^k)-\alpha_i^k\Bigg),
\end{align}
where $I_\mathcal{X}(\tilde{x})$ denotes the indicator function 
\begin{align}
I_\mathcal{X}(\tilde{x})=\begin{cases} 0,&\tilde{x}\in\mathcal{X},\nn\\ +\infty,&\tilde{x}\notin\mathcal{X}.\end{cases}
\end{align}

Comparing with the CADMM update in (\ref{eqn:distributedversion}), we see that (\ref{eqn:ConstrainedDOpt}) (and hence (\ref{eqn:ProjCADMMa})) together with (\ref{eqn:ProjCADMMb}) is the CADMM update with $\mathcal{T_X}(x_i^k)$ and $\alpha_i^k$ as local variables and $\frac{1}{2}(\tilde{x}-r_i)^2+I_\mathcal{X}(\tilde{x})$ as the local objective function at node $i$. As such, the convergence of the original CADMM implies that (\ref{eqn:ProjCADMMa}) and (\ref{eqn:ProjCADMMb}) yield 
\begin{align}
\lim_{k\to\infty}\mathcal{T_X}(x_i^k)&=\argmin_{\tilde{x}}\sum_{i=1}^n\left(\frac{1}{2}(\tilde{x}-r_i)^2+I_\mathcal{X}(\tilde{x})\right)\nn\\&=\mathcal{T_X}(\bar{r}).\nn
\end{align}
Thus, BQ-CADMM can also be viewed from applying the usual rounding quantizer to the original CADMM update of a constrained least-squares problem.

To recap, BQ-CADMM can be interpreted as  applying bounded quantization to local variables of an unconstrained least-squares problem, or alternatively, applying unbounded quantization to updates of a constrained least-squares problem. If $L$ is chosen to be large enough such that $\mathcal{T_X}(x_i^k)=x_i^k$ for all $k$, then $\mathcal Q_b(x_i^k)=\mathcal Q(x_i^k)$ and BQ-CADMM becomes the same as our previous algorithm DQ-CADMM in \cite{Zhu2016TSP,Zhu2016TSPerr}. In this sense, DQ-CADMM can be regarded as a special case of BQ-CADMM. The introduced projection operator or the additional indicator function, however, makes BQ-CADMM much more complicated to analyze. In \cite{Zhu2009,Zhu2016TSP, Zhu2016TSPerr} where dithered quantizer is studied, convergence comes from the linearity of the CADMM update for average consensus (cf. Eq.~(\ref{eqn:iterateform})) as well as the fact that the expectation of the output of a dithered quantizer is equal to the input. This approach fails to apply to BQ-CADMM because $\mathcal{T_X}(\cdot)$ is not a linear operator in general and one can not simply change the order of  projection and expectation operations. Additionally,  (\ref{eqn:ProjCADMMa}) and (\ref{eqn:ProjCADMMb}) no longer possess the linear convergence rate due to the introduced bounded constraint. Thus, the idea of using the linear convergence rate, as in \cite{Zhu2016TSP,Zhu2016TSPerr}, does not work for BQ-CADMM. On the other hand, if we view BQ-CADMM as a result from the CADMM update on the unconstrained least-squares problem modified by the bounded quantizer $\mathcal{Q}_b(\cdot)$, it is still unclear how BQ-CADMM performs since the quantization error of $\mathcal{Q}_b(\cdot)$ can be unbounded for a real input. Fortunately, local variables $\alpha_i^k$'s are inherently bounded by the BQ-CADMM update, as stated in the following lemma.
\begin{lemma}
\label{lem:alphabd}
Given local data $r_i\in\mathbb{R}$, consider BQ-CADMM with $x_i^0=0$ and $\alpha_i^0=0$ for $i=1,2,\ldots,n$. Then $\alpha_i^k$ is finitely bounded by
\begin{align}
|\alpha_{i}^k|\leq(1+6\rho|\mathcal{N}_i|)L+|r_i|.\nn 
\end{align}
\end{lemma}
\begin{IEEEproof}
See Appendix.
\end{IEEEproof}

To establish convergence result of BQ-CADMM, we define $\bm s_Q^k=[\mathcal{Q}_b(\bm x^k);\bm\alpha^k;\bm r]$ and rewrite BQ-CADMM as the standard CADMM update on $\bm s_Q^k$ plus an error term caused by bounded quantization. We first have 
\begin{align}
\mathcal{Q}_b(\bm x^{k+1})&= \bm x^{k+1} + \bm e_b(\bm x^{k+1}),\nn\\ 
\bm\alpha^{k+1}&=\bm\alpha^k+\rho\bm L_-\bm x^{k+1}+\rho\bm L_-\bm e_b(\bm x^{k+1}),\nn
\end{align}
where $\bm e_b(\bm x^{k+1})=\mathcal{Q}_b(\bm x^{k+1})-\bm x^{k+1}$ denotes the quantization error. By defining $\bm s_e^{k+1}=[\bm e_b(\bm x^{k+1});\rho\bm L_-\bm e_b(\bm x^{k+1});\bm 0]$, we have the BQ-CADMM update equivalent to
\begin{align}
\label{eqn:sstandarderr}
\bm s_Q^{k+1}=\bm D\bm s_Q^k+\bm s_e^{k+1}\triangleq\Phi(\bm s_Q^k).
\end{align}
It is important to note that the above update $\Phi(\cdot)$ is deterministic, i.e., given $\bm s_1=\bm s_2$, we must have $\Phi(\bm s_1)=\Phi(\bm s_2)$. We will use this fact to establish the main result given below.
\begin{theorem}
\label{thm:mainresults}
For BQ-CADMM in Algorithm \ref{tab:QCADMM}, there exists a finite iteration $k_0>0$ such that for $k\geq k_0$ all the quantized variable values
\begin{itemize}
\item either converge to the same quantization value:
$$\mathcal{Q}_b(x_1^k)=\cdots=\mathcal{Q}_b(x_n^k)\triangleq x_Q^*,$$
where $x_{Q}^*\in\Lambda_\mathcal{X}$ and
\begin{align}
\label{eqn:ProjCE}
\left|x^*_{Q}-\mathcal{T_X}(\bar{r})\right|\leq \left(1+4\rho\frac{m}{n}\right)\frac{\Delta}{2}.
\end{align}

\item or cycle around the average $\bar{r}$ with a finite period $T\geq2$, i.e., $x_i^k=x_i^{k+T}$, $i=1,2,\ldots, n$. Furthermore, the sample average over one period reaches a consensus:
\begin{align}
\label{eqn:cyclic_consensus}
\frac{1}{T}\hspace{-1pt}\sum_{l=1}^{T}\hspace{-2pt}\mathcal{Q}_b\left(x_1^{k+l}\right)=\cdots =\frac{1}{T}\hspace{-1pt}\sum_{l=1}^{T}\hspace{-2pt}\mathcal{Q}_b\left(x_n^{k+l}\right)\triangleq\bar{x}_Q^{*},
\end{align}
and 
\begin{align}
\label{eqn:cyclic_bd}
\left|\bar{x}_Q^{*}-\bar{r}\right|\leq\left(1+4\rho\frac{m}{n}\right)\Gamma_0,
\end{align}
where 
\begin{align}
\label{eqn:Tau0}
\Gamma_0\triangleq\max\left\{\frac{\Delta}{2},\frac{4\rho n L}{1+2\rho n}\right\}.
\end{align}
\end{itemize}
\end{theorem}
\begin{IEEEproof}
We first show that the sequence $\bm s_Q^k$ is either convergent or cyclic. Note first that $\mathcal{Q}_b(x_i^k)$ is bounded by $L$ and can only be a multiple of $\Delta$. Similarly, $\alpha_i^k$ is bounded as per Lemma~\ref{lem:alphabd} and is a multiple of $\rho\Delta$ as seen from the $\alpha_i$-update of BQ-CADMM. Thus, the number of states of local variables $[\mathcal{Q}_b(x_i^k);\alpha_i^k]$ at node $i$ is bounded by 
\begin{align}
\label{eqn:statesBD}
&\hspace{-.5in}\left(\frac{2L}{\Delta}+1\right)\left(\frac{L+|r_i|}{\rho\Delta}+\frac{6|\mathcal{N}_i|L}{\Delta}\right)\nn\\
~~\leq&\left(\frac{2L}{\Delta}+1\right)\left(\frac{L+\max_i|r_i|}{\rho\Delta}+\frac{6nL}{\Delta}\right)\triangleq B.
\end{align}
Therefore, the number of possible states of $\bm s_Q^k$ is bounded by $B^n$ which is finite for given $r_i$'s because $\rho$ and $\Delta$ are both positive and fixed. From~(\ref{eqn:sstandarderr}), the BQ-CADMM update $\bm s_Q^{k+1}=\Phi (\bm s_Q^k)$ is a deterministic function of only $\bm s_Q^k$. Therefore, $\bm s_Q^k$ must be either convergent or cyclic with a finite period $T\geq 2$ after a finite number of iterations denoted as $k_0$. The rest of the proof is to establish the respective error bounds.

{\it Convergent case:~}In this case we know $\bm\alpha^{k+1}=\bm\alpha^k$ for $k\geq k_0$. Then the $\bm\alpha$-update of BQ-CADMM together with Lemma~\ref{lem:graphLemma} implies $$ \mathcal{Q}_b(\bm x^k)\in\mathcal{C}(\bm1_n).$$ Thus, $\mathcal{Q}_b(\bm x^k)$ reaches a consensus at a common quantization point in $\Lambda_\mathcal{X}$. Denote $x_i^k=x_i^*$ and $\alpha_i^k=\alpha_i^*$ for $k\geq k_0$, i.e., when convergence is reached. Further letting $e_i^*=x_Q^*-\mathcal{T_X}(x_{i}^*)$, we have $|e_i^*|\leq\frac{\Delta}{2}$ since $\mathcal{Q}_b(\cdot)=\mathcal{Q}\circ\mathcal{T_X}(\cdot)$ with $\mathcal{Q}(\cdot)$ being the usual rounding quantization. Taking $k\to\infty$ on both sides of (\ref{eqn:ConstrainedDOpt}) and using the optimality condition for minimizing a convex function (see, e.g., \cite{BoydConvexOpt}), we get
\begin{align}
\mathcal{T_X}(x_{i}^*)-r_i&+2\rho|\mathcal{N}_i|\mathcal{T_X}(x_{i}^*)+\partial I_\mathcal{X}\left(\mathcal{T_X}(x_{i}^*)\right)\nn\\
&~~~~~~~~-\Bigg(\rho|\mathcal{N}_i|x_{Q}^*
+\rho\sum_{j\in\mathcal{N}_i}x_{Q}^*-\alpha_{i}^*\Bigg)=0,\nn
\end{align}
where $\partial I_\mathcal{X}(\mathcal{T_X}(x_{i}^*))$ is a subgradient of $I_\mathcal{X}(\cdot)$ at $\mathcal{T_X}(x_{i}^*)$. After rearranging and plugging in $\mathcal{T_X}(x_{i}^*)=x_Q^*-e_i^*$, we obtain 
\begin{align}
\label{eqn:OptCon}
x_Q^*-e_i^*-r_i+\partial I_\mathcal{X}\left(x_Q^*-e_i^*\right)+\alpha_{i}^*=2\sum_{j\in\mathcal{N}_i}|\mathcal{N}_i|e_j^*.
\end{align}
Summing up both sides of (\ref{eqn:OptCon}) from $i=1$ to $n$ yields
\begin{align}
\label{eqn:sumOptCon}
\sum_{i=1}^n\left(x_Q^*-e_i^*-r_i+\partial I_\mathcal{X}\left(x_Q^*-e_i^*\right)\right)=2\sum_{i=1}^n|\mathcal{N}_i|e_i^*,
\end{align}
where $\sum_{i=1}^n\alpha_{i}^*=0$ is due to $[\alpha_1^*;\alpha_2^*;\cdots;\alpha^*_n]\in\mathcal{C}(\bm L_-)$ and Lemma \ref{lem:graphLemma}. Since $\mathcal{T_X}(\bar{r})=\argmin_{\tilde{x}}\sum_{i=1}^n \big(\frac{1}{2}(\tilde{x}-r_i)^2+I_\mathcal{X}(\tilde{x})\big)$, we also have\footnote{Here $\partial I_\mathcal{X}\left(\mathcal{T_X}(\bar{r})\right)$ denotes a subgradient, i.e., an element from the subdifferential of $I_\mathcal{X}(\cdot)$ at $\mathcal{T_X}(\bar{r})$; we simply use the same notation for all the subgradients despite the fact that they might be different values.}
\begin{align}
\label{eqn:OptConAvg}
\sum_{i=1}^n\left( \mathcal{T_X}(\bar{r})-r_i+\partial I_\mathcal{X}\left(\mathcal{T_X}(\bar{r})\right)\right)=0.
\end{align}
Subtracting (\ref{eqn:OptConAvg}) from (\ref{eqn:sumOptCon}), we get that
\begin{align}
\label{eqn:consensuserr}
&~\sum_{i=1}^n\left(\mathcal{T_X}(\bar{r})-x_Q^*+\partial I_\mathcal{X}\left(\mathcal{T_X}(\bar{r})\right)-\partial I_\mathcal{X}\left(x_Q^*-e_i^*\right)\right)\nn\\
=&~\sum_{i=1}^n\left(2|\mathcal{N}_i|+1\right)e_i^*.
\end{align}
We only consider the case where $|\mathcal{T_X}(\bar{r})-x_Q^*|>\frac{1}{2}\Delta$; otherwise (\ref{eqn:ProjCE}) holds trivially. Recall that $|e_i^*|\leq\frac{\Delta}{2}$, we get  $$(\mathcal{T_X}(\bar{r})-x_Q^*)\left(\mathcal{T_X}(\bar{r})-(x_Q^*-e_i^*)\right)> 0.$$ Note that the convexity of $I_\mathcal{X}(\cdot)$ also implies $$\left(\partial I_\mathcal{X}\hspace{-1pt}\left(\mathcal{T_X}(\bar{r})\right)-\partial I_\mathcal{X}\hspace{-1pt}\left(x_Q^*-e_i^*\right)\right)\hspace{-1pt}\left(\mathcal{T_X}(\bar{r})-(x_Q^*-e_i^*)\right)\geq 0.$$ 
The above two inequalities indicate that either $\partial I_\mathcal{X}\left(\mathcal{T_X}(\bar{r})\right)-\partial I_\mathcal{X}\left(x_Q^*-e_i^*\right)=0$, or $\partial I_\mathcal{X}\left(\mathcal{T_X}(\bar{r})\right)-\partial I_\mathcal{X}(x_Q^*-e_i^*)$ and $\mathcal{T_X}(\bar{r})-x_Q^*$ have the same sign. Then the following is true:
\begin{align}
&\left|\mathcal{T_X}(\bar{r})-x_Q^*+\partial I_\mathcal{X}\left(\mathcal{T_X}(\bar{r})\right)-\partial I_\mathcal{X}\left(x_Q^*-e_i^*\right)\right|
\nn\\
\geq&\left|\mathcal{T_X}(\bar{r})-x_Q^*\right|\nn.
\end{align}
Together with (\ref{eqn:consensuserr}), we can establish the upper bound
\begin{align}
\left|x_Q^*-\mathcal{T_X}(\bar{r})\right|
\leq&~\frac{1}{n}\left|\sum_{i=1}^n\left(2\rho|\mathcal{N}_i|+1\right)e_i^*\right|\nn\\
\leq&~\left(1+\rho\frac{4m}{n}\right)\frac\Delta2\nn,
\end{align} 
where we use the fact that $\sum_{i=1}^n|\mathcal{N}_i|=2m$ for an undirected connected graph.

{\it Cyclic case:~}That $\bm s_Q^k$ cycles with a period $T$ implies $$\bm\alpha^{k+T}-\bm\alpha^k=\rho\bm L_-\sum_{l=1}^T\mathcal{Q}_b(\bm x^{k+l})=\bm 0,$$
which then leads to (\ref{eqn:cyclic_consensus}) by Lemma~\ref{lem:graphLemma}. While the bound on $\alpha_i^k$ in Lemma~\ref{lem:alphabd} imposes a bound on $x_i^k$ through the $x_i$-update, the cyclic behavior itself can be used to derive a tighter error bound.  

Consider the local variable $x_i^k$ over one period. When $x_i^k\in\mathcal{X}$ for the entire period (and hence for all $k\geq k_0$), we simply have $|e_b(x_i^k)|\leq\frac{\Delta}{2}$. If there is some $x_i^{k+1}>L$ with $k\geq k_0+T-1$, there must exist a $k'\leq k$ which is the largest index such that $\mathcal{Q}_b(x_i^{k'})<L$ and $\mathcal{Q}_b(x_i^{k'+1})=L$, for otherwise $\mathcal{Q}_b(x_i^k)$ converges and hence BQ-CADMM must converge due to (\ref{eqn:cyclic_consensus}). Then we have $\mathcal{Q}_b(x_i^{k'+l})=L $ for $l=1,2,\ldots,k+1-k'$. Recall the $\alpha_i$-update, we get $\alpha_i^{k'+l}\leq\alpha_i^{k'}$ as $L$ is the largest quantization value. We further write the $x_i^{k+1}$-update as
\begin{align}
\label{eqn:xkbdcyc}
&\hspace{2pt}x_i^{k+1}\nn\\
=&\hspace{2pt}\frac{1}{1+2\rho|\mathcal{N}_i|}\Bigg(\rho|\mathcal{N}_i|\mathcal{Q}_b(x_i^k)+\rho\sum_{j\in\mathcal{N}_i}\mathcal{Q}_b(x_j^k) -\alpha_i^k + r_i\Bigg)\nn\\
\stackrel{(a)}{\leq}&\hspace{2pt}\frac{1}{1+2\rho|\mathcal{N}_i|}\left(2\rho|\mathcal{N}_i|L-\alpha_i^{k'} + r_i \right)\nn\\
\stackrel{(b)}{=}&\hspace{2pt}x_i^{k'}+\frac{\rho}{1+2\rho|\mathcal{N}_i|}\Bigg(2|\mathcal{N}_i|L+\sum_{j\in\mathcal{N}_i}\mathcal{Q}_b(x_j^{k'})-|\mathcal{N}_i|\mathcal{Q}_b(x_i^{k'}\hspace{-1pt})\nn\\
&\hspace{2pt}-|\mathcal{N}_i|\mathcal{Q}_b(x_i^{k'-1})-\sum_{j\in\mathcal{N}_i}\mathcal{Q}_b(x_j^{k'-1})\Bigg)\nn\\
\stackrel{(c)}{\leq}&\hspace{2pt}\frac{5\rho|\mathcal{N}_i|L}{1+2\rho|\mathcal{N}_i|}+x_i^{k'}-\frac{\rho|\mathcal{N}_i|}{1+2\rho|\mathcal{N}_i|}\mathcal{Q}_b(x_i^{k'}),
\end{align}
where $(a)$ and $(c)$ are due to the fact that $-L\leq\mathcal{Q}_b(\tilde{x})\leq L$ for any $\tilde{x}\in\mathbb{R}$, and $(b)$ is from the BQ-CADMM update at the $k'$th iteration. Since $k'$ is the index such that $\mathcal{Q}_b(x_i^{k'})<L$, it is straightforward to see that (\ref{eqn:xkbdcyc}) takes the largest value at $x_i^{k'}=L-\frac{\Delta}{2}$. This implies for $k\geq k_0$,
\begin{align}
x_i^k&\leq L-\frac{\Delta}{2}+\frac{5\rho|\mathcal{N}_i|L}{1+2\rho|\mathcal{N}_i|}-\frac{\rho|\mathcal{N}_i|(L-\Delta)}{1+2\rho|\mathcal{N}_i|}\nn\\
&< L+\frac{4\rho|\mathcal{N}_i|L}{1+2\rho|\mathcal{N}_i|}.\nn
\end{align}
Similarly, it can be shown that for $k\geq k_0$, $$x_i^k>-L-\frac{4\rho|\mathcal{N}_i|L}{1+2\rho|\mathcal{N}_i|}.$$
Therefore, when BQ-CADMM reaches a cyclic state, we have
\begin{align}
\label{eqn:xbdcyclic2}
|x_i^k|\leq L+\frac{4\rho|\mathcal{N}_i|L}{1+2\rho|\mathcal{N}_i|},
\end{align}
and the quantization error satisfies
\begin{align}
\label{eqn:errbd_cyc}
\left|e_b(x_i^k)\right|\leq\max\left\{\frac{\Delta}{2}, \frac{4\rho|\mathcal{N}_i|L}{1+2\rho|\mathcal{N}_i|}\right\}\leq \Gamma_0,
\end{align}
where $\Gamma_0$ is defined in (\ref{eqn:Tau0}) and we use the fact that $1\leq |\mathcal{N}_i|< n$.

To derive the error bound, we now summarize the local variable values over one period and use (\ref{eqn:errbd_cyc}), which leads to 
\begin{align}
\label{eqn:ebdcyc}
\left|\frac{1}{T}\sum_{l=1}^{T}\mathcal{Q}_b(x_i^{k+l})-\frac{1}{T}\sum_{l=1}^{T}x_i^{k+l}\right|= \left|\bar{x}_Q^*-\frac{1}{T}\sum_{l=1}^{T}x_i^{k+l}\right|
\leq \Gamma_0.
\end{align}
By the $x_i$-update, we can also get
\begin{align}
\label{eqn:sumsumsum}
&(1+2\rho|\mathcal{N}_i|)\frac{1}{T}\sum_{l=1}^{T}x_i^{k+l}-\rho|\mathcal{N}_i|\frac{1}{T}\sum_{l=1}^{T}\mathcal{Q}_b(x_i^{k+l})\nn\\
&~~-\rho\sum_{j\in\mathcal{N}_i}\left(\frac{1}{T}\sum_{l=1}^{T}\mathcal{Q}_b(x_i^{k+l})\right)+\sum_{l=1}^{T}\alpha_i^{k+l}-r_i=0.
\end{align}
Summing both sides of (\ref{eqn:sumsumsum}) over $i$ and using (\ref{eqn:cyclic_consensus}), we get 
\begin{align}
\label{eqn:just11}
&\sum_{i=1}^nr_i+\sum_{i=1}^n\left((1+2\rho|\mathcal{N}_i|)\frac{1}{T}\sum_{l=1}^{T}x_i^{k+l}\right)\nn\\
&~~~~~~~~~~~~~~~~-\sum_{i=1}^n\left(2\rho|\mathcal{N}_i|\frac{1}{T}\sum_{l=1}^{T}\mathcal{Q}_b(x_i^{k+l})\right)=0.
\end{align}
Finally, using (\ref{eqn:ebdcyc}) and dividing both sides of (\ref{eqn:just11}) by $n$, we can bound the consensus error by
\begin{align}
\left|\bar{x}_Q^*-\bar{r}\right|\leq\left(1+4\rho\frac{m}{n}\right)\Gamma_0.\nn
\end{align}
This completes the proof.
\end{IEEEproof}

\begin{remark}
\label{rmk:Qitself}
We have mentioned that BQ-CADMM uses $\mathcal{Q}_b(x^k_i)$ for the $({k+1})$th update at agent $i$ even though agents can compute and store real values with infinite precision. The reason is to guarantee that $\bm\alpha^k$ is bounded and also lies in $\mathcal{C}(\bm L_-)$ in order to achieve a consensus result. Theorem \ref{thm:mainresults} may fail if $\bm\alpha^k\notin\mathcal{C}(\bm L_-)$.
\end{remark}
\begin{remark}
\label{rmk:InitialState}
While BQ-CADMM in Algorithm~\ref{tab:QCADMM} initializes $\bm x^0=\bm\alpha^0=\bm0$, one can use similar arguments to show the same convergence result of Theorem~\ref{thm:mainresults} for any $\bm x^0\in\mathbb{R}^n$ and $\bm\alpha^0\in\mathcal{C}(\bm L_-)$. As a result, we can pick $\bm x^0$ and $\bm \alpha^0$ that are respectively closer to the optima $\bm1_n\bar{r}$ and $\bm r- \bm1_n\bar{r}$, which usually leads to better consensus performance. Since Theorem~\ref{thm:mainresults} holds for any $\rho>0$, we can also use decreasing algorithm parameter $\rho$ such that the algorithm proceeds fast in early stages and guarantees certain consensus accuracy upon convergence or cycling (cf.~Section~\ref{sec:DSforPS}).
\end{remark}
\begin{remark}
\label{rmk:ceBQ} 
With convergence or small enough $\rho$ such that $\Gamma_0\leq\frac{\Delta}{2}$, the error bounds (\ref{eqn:ProjCE}) and (\ref{eqn:cyclic_bd}) are exactly the same as those of DQ-CADMM when $\bar{r}\in\mathcal{X}$. Note that even with the prior knowledge that $\bar{r}\in\mathcal{X}$, the algorithms of \cite{Kar2010,Li2011} fail to provide a guaranteed consensus result as agents' data can still be unbounded.  This not only saves energy for data communication but also broadens its applications; see consensus based detection in \cite{Zhu2016Detection} where we show that BQ-CADMM with one-bit information exchange between linked nodes at each iteration suffices to achieve the optimal exponential decay of error probabilities under the maximum {\it a posterior} and the Neyman-Pearson criteria. 
\end{remark}

While the proposed algorithm is not guaranteed to converge for all cases, we can use history of the agent variable values (e.g., the running average technique) to ensure asymptotic convergence at a consensus. As such, we refer to $x_Q^*$ in the convergent case and $\bar{x}_Q^*$ in the cyclic case as the resulting consensus value. As a direct result from Theorem~\ref{thm:mainresults}, the following corollary states that BQ-CADMM must converge when the data average resides outside the bounded set $\mathcal{X}$.
\begin{corollary} BQ-CADMM must converge to a consensus at a quantization level in $\Lambda_\mathcal{X}$ when 
$$|\bar{r}|-L>\left (1+4\rho\frac{m}{n}\right)\Gamma_0.$$
If we further pick $\rho<\frac{n}{4m}$, then $\mathcal{Q}_b(x_i^k)$ must converge to $\text{sgn}(\bar{r})L$ where $\text{sgn}(\bar{r})=1$ if $\bar{r}>0$ and $\text{sgn}(\bar{r})=-1$ if $\bar{r}<0$.

\end{corollary}

\subsection{Extended BQ-CADMM}
The error bounds in Theorem \ref{thm:mainresults} indicate that BQ-CADMM eventually reaches a neighborhood of $\bar{r}$. When $\bar{r}$ is far from the bounded set $\mathcal{X}$, however, we may suffer from a large consensus error. Pick $\rho$ small enough such that $\Gamma_0\leq\frac{\Delta}{2}$ and assume that $L\gg \left({1}+4\rho\frac{m}{n}\right)\frac{\Delta}{2}$ (e.g, $L\geq 5\left({1}+4\rho\frac{m}{n}\right)\frac{\Delta}{2}$). Then Corollary~1 implies that $\bar{r}$ is either outside $\mathcal{X}$ or around its boundary when BQ-CADMM converges at $|x_Q^*|=L$. If each agent subtracts $\text{sgn}(x_Q^*)L$ from $r_i$, the new data average becomes $|\bar{r}-\text{sgn}(x_Q^*)L|=\left||\bar{r}|-L\right|<|\bar{r}|$. Thus, we can run BQ-CADMM with updated data $r_i-\text{sgn}(x_Q^*)L$. This process is repeated until BQ-CADMM either converges to $|x_Q^*|\neq L$ or reaches a cyclic result. The consensus value at each agent is simply the sum of the final consensus value and the offset from each BQ-CAMM call. The above procedure is referred to as extended BQ-CADMM (EBQ-CADMM) presented in Algorithm~\ref{tab:EBQCADMM}.
\begin{Algorithm}{\linewidth}
	\caption{EBQ-CADMM for distributed average consensus}
	\begin{algorithmic}[1]
	\label{tab:EBQCADMM}
	\REQUIRE Initialize~$t_i=0$ for each agent $i,i=1,2,\ldots,n$. Pick $\rho>0$ such that $\Gamma_0\leq\frac{\Delta}{2}$.
	\WHILE{true}
		\STATE Run BQ-CADMM with data $r_i$ at agent $i$ such that either convergence or cycling is achieved. Denote the consensus value as $x_{BQ}$.
		\IF {$|x_{BQ}|=L$}
		\STATE set $t_i=t_i+\text{sgn}(x_{BQ})L$ and $r_i=r_i-\text{sgn}(x_{BQ})L$	
		\ELSE
		\STATE break
		\ENDIF	
	\ENDWHILE
	\RETURN $t_i$ and $x_{BQ}$
	\end{algorithmic}
\end{Algorithm}

We remark that an additional BQ-CADMM is run only when a convergence is reached at either $-L$ or $L$. Notice that consensus is always guaranteed by Theorem~\ref{thm:mainresults}; that is, if a node converges at $-L$ (or $L$), then every other node converges at $-L$ (or $L$). Hence, global synchronization is not needed at each BQ-CADMM call. Since  EBQ-CADMM calls the BQ-CADMM at most $\left\lceil{\left|\bar{r}\right|}/{L}\right\rceil+1$ times, we have the following theorem directly from Theorem \ref{thm:mainresults} and Corollary~1.
\begin{theorem}
\label{thm:EBQCADMM}
Assume that $L\gg\left(1+4\rho\frac{m}{n}\right)\frac{\Delta}{2}$. For any agents' data $r_i$'s, EBQ-CADMM in Algorithm~\ref{tab:EBQCADMM} yields that $$t_1=\cdots=t_n\triangleq t^*,$$  and
\begin{align}
\left|t^*+x_{BQ}-\bar{r}\right|\leq\left({1}+4\rho\frac{m}{n}\right)\frac{\Delta}{2}.\nn
\end{align}
\end{theorem}

It is noted that the above error bound in Theorem~\ref{thm:EBQCADMM} has a non-vanishing term $\frac{\Delta}{2}$, which implies that the consensus result can not be made arbitrarily close to the true average in general. Nevertheless, by employing deterministic uniform quantization and a fixed step size, EBQ-CADMM possesses local variables which can be used to recover the original data average. Denote the final iteration index by $k'$. At the end of EBQ-CADMM, the last BQ-CADMM run either converges at some quantization point in $(-L,L)$, which indicates $|x_i^{k}|<L$, or cycles with $|x_i^{k}|\leq L+\frac{4\rho|\mathcal{N}_i|L}{1+2\rho|\mathcal{N}_i|}$ (cf.~Eq.~(\ref{eqn:xbdcyclic2})) for any $k\geq k'$. Thus, at each node $i$, we have 
\begin{align}
&\hspace{2pt}\left|-\alpha_i^{k'}+r_i-t^*\right|\nn\\
=&\hspace{2pt}\left|(1+2\rho|\mathcal{N}_i|)x_i^{{k'}+1} + \rho|\mathcal{N}_i|\mathcal{Q}_\delta(x_i^{k'})+\rho\sum_{j\in\mathcal{N}_i}\mathcal{Q}_\delta(x_i^{k'})\right|\nn\\
\leq&\hspace{2pt}L +8\rho|\mathcal{N}_i|L.\nn
\end{align}
Further picking $\rho=\frac{1}{8n}$, we can get $|-\alpha_i^{k'}+r_i-t^*|\leq2L$.
Also note that $\sum_{i=1}^n\alpha_i^{k'}=0$ according to Lemma~\ref{lem:graphLemma} and that $t^*$ is locally known to each node. Thus, setting local data to $-\alpha_i^{k'}+r_i-t^*$, existing quantized consensus algorithms can use finite-bit quantization to asymptotically achieve the exact average consensus; more precisely, the variable value at each node is close to $$\frac{1}{n}\sum_{i=1}^n\left(-\alpha_i^{k'}+r_i-t^*\right)=\bar{r}-t^*,$$
and the true average $\bar{r}$ is asymptotically achieved by adding~$t^*$.
\section{Algorithm Parameter and Stopping Criterion}
\label{sec:rho}
In this section, we discuss the effect of the algorithm parameter $\rho$ on BQ-CADMM and EBQ-CADMM. Since $L$ can be large enough such that $\mathcal{Q}_b(x_i^k)=\mathcal{Q}(x_i^k)$, the discussion also covers DQ-CADMM using the usual rounding quantizer in \cite{Zhu2016TSP, Zhu2016TSPerr} where the effect of $\rho$ is not explored.
\subsection{Consensus Error}
As seen from the consensus error bounds in (\ref{eqn:ProjCE}) and (\ref{eqn:cyclic_bd}), it is clear that they increase in the algorithm parameter $\rho$. With the knowledge of the number of nodes or its upper bound, one can pick $\rho$ small enough to achieve desired consensus accuracy. For example, if $\rho\leq\min\left\{\frac{n}{4m},\frac{\Delta}{8nL}\right\}=\frac{\Delta}{8nL}$ as $2m\leq{n(n-1)}$ and $L\geq\Delta$, the resulting consensus value of EBQ-CADMM is within one quantization resolution of the desired average. We also remark that practical consensus error is usually much smaller than the error bound (see simulations). 
\subsection{Cyclic Period}
\label{sec:CP}
In the proof of Theorem~\ref{thm:mainresults}, we show that the number of possible states of $\bm s_Q^k$ is bounded by $B^n$ where $B$ is given in (\ref{eqn:statesBD}). This bound works for both convergence time and cyclic period but is generally too loose to use in practice. Our simulations in Section~\ref{sec:CP} show that BQ-CADMM converges in most cases, particularly with small enough $\rho$. When the algorithm indeed cycles, the period depends on the network structure, agents' data as well as the algorithm parameter, with a range between $2$ and $15$ for all simulated cases. We also observe that the cyclic period in all simulations consists of two consecutive quantization levels for each node. Indeed, we can derive tighter consensus error bounds if the cyclic period only has two consecutive quantization levels, e.g., we can replace $L$ with $\Delta$ in $\Gamma_0$ that is defined by (\ref{eqn:Tau0}).

\subsection{Convergence Time}
\label{sec:rhoCT}
We refer to the smallest $k_0$ in convergent cases or smallest $k_0+T$ in cyclic cases as the convergence time for BQ-CADMM. First notice that it is meaningless to consider the convergence time without any constraint on the consensus error. To see this, recall that BQ-CADMM has initial variable values $x_i^0=\alpha_i^0=0$ at each node $i$ and $x_i^{1}=\frac{1}{1+2\rho|\mathcal{N}_i|} r_i$ at the first iteration. If $\rho$ is chosen large enough, e.g., $\rho>\max_i\frac{|r_i|}{\Delta}$, such that $|x_i^{1}|<\frac{\Delta}{2}$, then we have $\mathcal{Q}_b(x_i^k)=\alpha_i^k=0, k=0,1,\ldots,$ for all $i$ and the convergence time is $k_0=1$. Thus, we consider accelerating BQ-CADMM under some consensus accuracy constraint.

A direct approach is to apply similar techniques that accelerate the original ADMM. For example, reference \cite{Shi2014} studies the optimal parameter selection of $\rho$ that maximizes the convergence rate $\delta$ given in Theorem~\ref{thm:linearconvergence} and reference \cite{Giselsson2016} uses preconditioning to improve the convergence rate. However, these approaches require the knowledge of the network structure, which might be unrealistic for large scale networks. Furthermore, the resulting parameter selection may not meet the consensus accuracy requirement and indeed does not provide the best practical performance. Here we propose a heuristic selection $\rho =\frac{n}{m}$ which only requires (an estimate of) the number of nodes and the number of edges. Our intuition is based on the fact that a larger $\frac{m}{n}$ indicates that an agent on the average has more neighboring information for its update. Therefore, a smaller $\rho$ provides adequate updates towards the consensus. Conversely, when agents have less information available, a larger $\rho$ can help accelerate the speed. The good performance of this selection is validated by simulations. When $\rho=\frac{n}{m}$ fails to satisfy a given consensus accuracy constraint, a decreasing strategy can be used where the initial value of $\rho$ is set at $\frac{n}{m}$ and decreases as the number of iterations increases. In the case where only $n$ is known but $m$ is not, the decreasing strategy may simply start with $\rho=1$. 
\subsection{Stopping Criterion}
A natural approach for stopping  BQ-CADMM is to set the maximum number of iterations at each node, which then requires characterization of the convergence time. While our current analysis does not have a closed-form result,  we alternatively provide a rough upper bound based on our observations. At the beginning of the algorithm when local agent variables are far from the optima, i.e., when $\|\bm u^k-\bm u^*\|_{\bm G}^2$ is large, the CADMM update in (\ref{eqn:ProjCADMMa}) and (\ref{eqn:ProjCADMMb}) will make much progress towards the optima and the effect of the quantization operation $\mathcal{Q}(\cdot)$ is negligible. Thus, BQ-CADMM and the original CADMM for the constrained least-squares problem have similar trajectories at early stages. When local variables becomes close to the optima, the original CADMM still has monotone decreasing $\|\bm u^k-\bm u^*\|_{\bm G}^2$ while BQ-CADMM either converges or cycles due to the quantization operation. This behavior is observed in all our simulations and an example is provided in Fig.~\ref{fig:Comp}. As such, a rough estimate can be obtained based on the convergence rate result of the original CADMM when applied to the constrained least-squares problem. Here we use the result from \cite{He2015non} where $\|\bm u^{k+1}-\bm u^k\|_{\bm G}^2$ is viewed as a metric of how close local variables are to the optima. The convergence rate of the original CADMM is 
\begin{align}
\label{eqn:CADMMrate}
\hspace{2pt}\|\bm u^{k+1}-\bm u^{k}\|_{\bm G}^2
\leq\hspace{2pt}\frac{1}{k+1}\|\bm u^0-\bm u^*\|_{\bm G}^2.
\end{align}
We can bound the right-hand side of (\ref{eqn:CADMMrate}) by a small value and our simulations show that $4\rho\Delta^2$ is a good choice. Thus, a rough bound on convergence time, with $x_i^0=\alpha_i^0=0$ for $i=1,2,\ldots,n$, is given by 
$$\frac{1}{\Delta^2}\left( n\mathcal{T_X}(\bar{r})^2+\frac{1}{2\rho^2{\lambda_2(\bm L_-)}}\|\bm r-\bm1_n\mathcal{T_X}(\bar{r})\|_2^2\right).$$
We remark that this rough bound is quite loose in most of our simulations; see Fig.~\ref{fig:CT}. An upper bound on the convergence time, which depends on the network topology and agents' data, may still be insufficient as these quantities are locally unknown. In fully distributed settings, we can run additional algorithms to determine if a consensus has been reached; see, e.g., \cite{Manitara2016,Yadav2007}. These additional algorithms may take a long running time in large networks and will require extra data communications. In general, finding a practically meaningful and efficient stopping criterion in a fully distributed manner remains open.

As to determining the cyclic behavior, nodes may record a certain number of consecutive variable values and check if any cycle exists. Though no bound on cyclic period is provided in this paper, our simulations show that the cyclic period is between $2$ and $15$. In summary, we can set the maximum number of iterations or run the algorithm as long as permitted. If there is no convergence which can be determined by additional algorithms, we can then run other algorithms with agents' data $-\alpha_i^k+r_i$ which are bounded with a known bound at each node.

\section{Simulations}
\label{sec:simulation}
This section investigates the proposed algorithms via numerical examples. We construct a connected network with $n$ nodes and $m$ edges by first generating a complete graph of $n$ nodes and then randomly removing $\frac{n(n-1)}{2}-m$ edges while ensuring the graph stays connected. Let $\Delta = 1$ throughout this section. 
\subsection{Consensus Error}
\label{sec:consensuserr}
\begin{figure}[htbp]
	\centering
	\includegraphics[width=\linewidth]{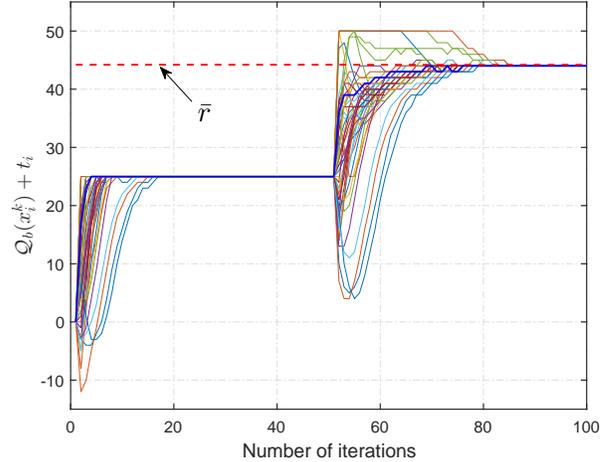}%
	\caption{Trajectories of EBQ-CADMM; $n=50$, $m=100$, $L=\frac{n}{2}$, $\Delta=1$, $\rho=0.5$, and $r_i\sim\mathcal(n,n^2)$.}
	\label{fig:EBQCADMM}
\end{figure}

We first illustrate how BQ-CADMM and EBQ-CADMM proceed by showing the trajectories of agents' variables. We simulate a connected network with $n=50$ nodes and $m=100$ edges. Let $L = \frac{n}{2}=25$, $\rho = \frac{n}{m}=0.5$, and $r_i\sim\mathcal{N}(n,n^2)$. Set also the maximum iteration number of each BQ-CADMM to $50$. By this setting, it is very likely to have $|\bar{r}|>L$ and we thus apply EBQ-CADMM for data averaging. The trajectories of $\mathcal{Q}_b(x_i^k)+t_i,i=1,2,\ldots,n$ are plotted in Fig.~\ref{fig:EBQCADMM}. In this example, the desired average is $44.20$ and the resulting consensus value of EBQ-CADMM is $44$. The consensus error is $0.20$, which is much smaller than the error bound $(1+4\rho\frac{m}{n})\frac{\Delta}{2}= 2.5$. The figure also indicates that EBQ-CADMM calls BQ-CADMM twice with the first and second calls converging at iterations $17$ and $33$, respectively. As a note, we run this setup for $10,000$ times and no cyclic result is observed.

\begin{figure}[htbp]
	\centering
		\includegraphics[width=\linewidth]{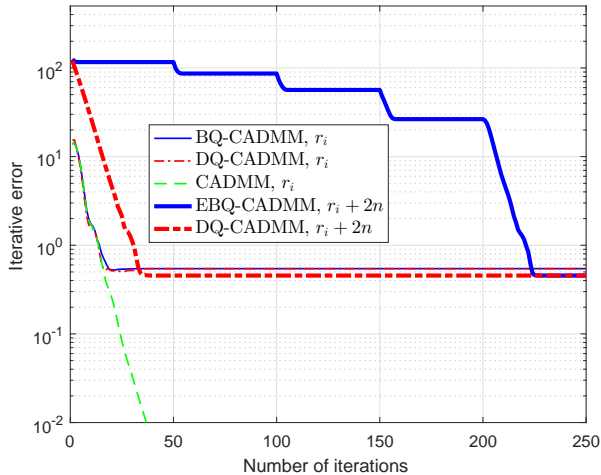}
	
	\caption{Iterative errors of BQ-CADMM, EBQ-CADMM, DQ-CADMM, and CADMM; $n=75$, $m=200$, $L=30$, $\Delta=1$, $\rho=0.5$, and $r_i\sim\mathcal{N}(0,n^2)$.}
	\label{fig:Comp}
\end{figure}

We next compare the proposed algorithms with DQ-CADMM that uses rounding quantizer $\mathcal{Q}(\cdot)$ to investigate the effect of the bounded constraint. Denote $\bm t=[t_1;t_2;\cdots;t_n]$ for EBQ-CADMM and define the iterative error as 
\begin{align}
\label{eqn:bounddQ}
\begin{cases}
{\left\|\mathcal{Q}_b({\bm x}^k)-\bm1_n\bar{r}\right\|_2}/{\sqrt{n}},&~\text{BQ-CADMM},\\
{\left\|\mathcal{Q}_b({\bm x}^k)+\bm t-\bm1_n\bar{r}\right\|_2}/{\sqrt{n}},&~\text{EBQ-CADMM},\\
{\left\|\mathcal{Q}({\bm x}^k)-\bm1_n\bar{r}\right\|_2}/{\sqrt{n}},&~\text{DQ-CADMM},\\
{\left\|\bm x^k - \bm1_n\bar{r}\right\|/\sqrt{n}},&~\text{CADMM},\nn\end{cases}
\end{align} which is equal to the consensus error when a consensus is reached. Set $n=75$, $m=200$ and $L =30$. Pick $\rho=0.5$ and set the maximum iteration number of each BQ-CADMM call in EBQ-CADMM to $50$. We consider two cases: one with the average in the bounded set and the other with the average outside the set. Specifically, we generate $r_i\sim\mathcal{N}(0,n^2)$ and run the distributed averaging algorithms twice with data $r_i$ and $r_i+2n$, respectively. We use BQ-CADMM for the former case since the average lies in the bounded set with high probability. Simulation result is presented in Fig.~\ref{fig:Comp}.

From Fig.~\ref{fig:Comp}, all quantized consensus algorithms converge to a consensus after certain iterations with consensus errors much smaller than the upper bound, computed to be $1.83$. When $\bar{r}\in\mathcal{X}$, BQ-CADMM and DQ-CADMM perform similarly and their iterative errors decrease exponentially until consensus is reached, while the original CADMM makes the iterative error decrease to $0$. When $|\bar{r}|>L$, DQ-CADMM converges much faster than EBQ-CADMM. This is because each BQ-CADMM call has iterations running even when a consensus is reached. 


\subsection{Cyclic Period}
\begin{figure*}[htbp]
	\centering
	\subfigure[]{%
		\includegraphics[width=0.325\linewidth,height=0.26\linewidth]{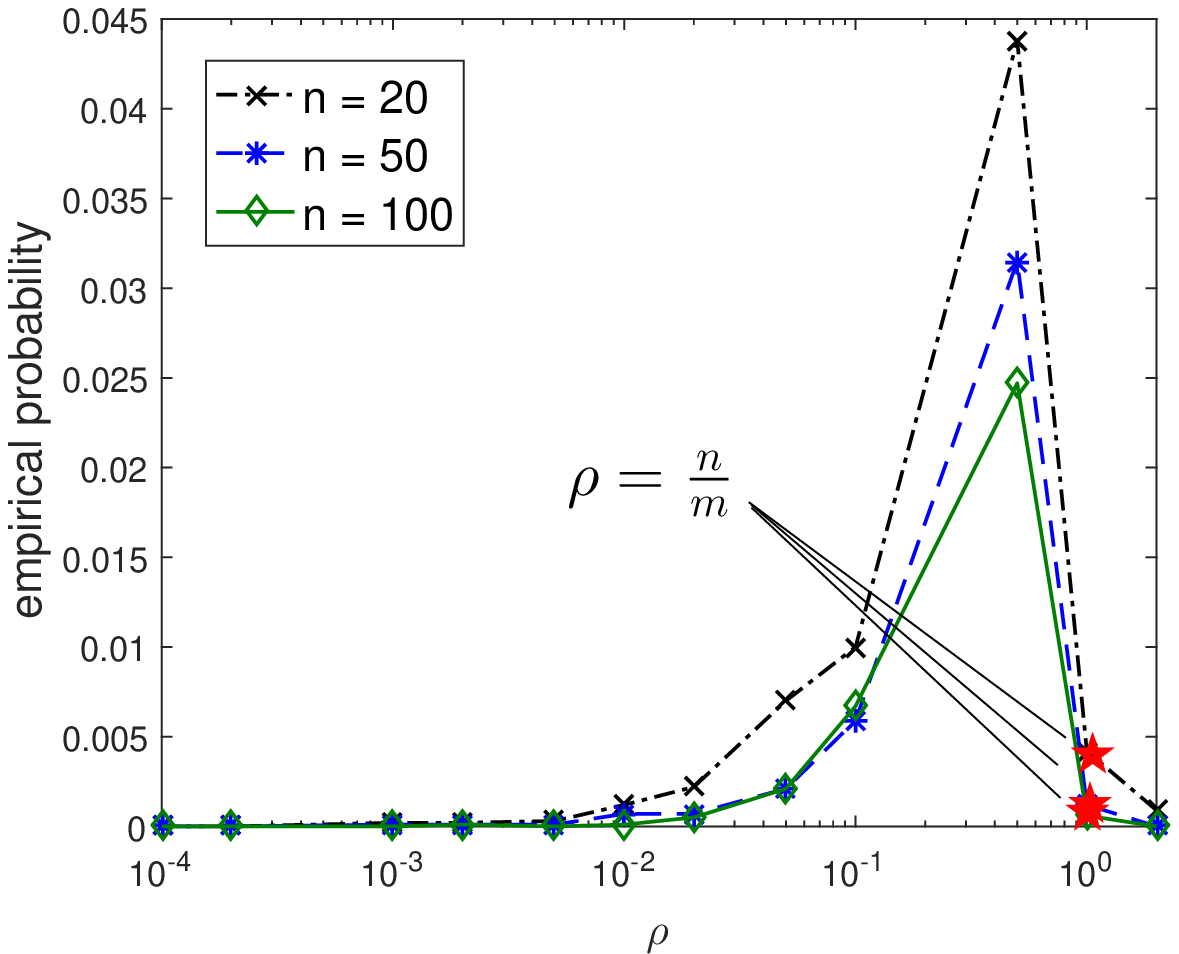}%
		\label{fig:Cyctar}%
	}
	\subfigure[]{%
		\includegraphics[width=0.325\linewidth,height=0.26\linewidth]{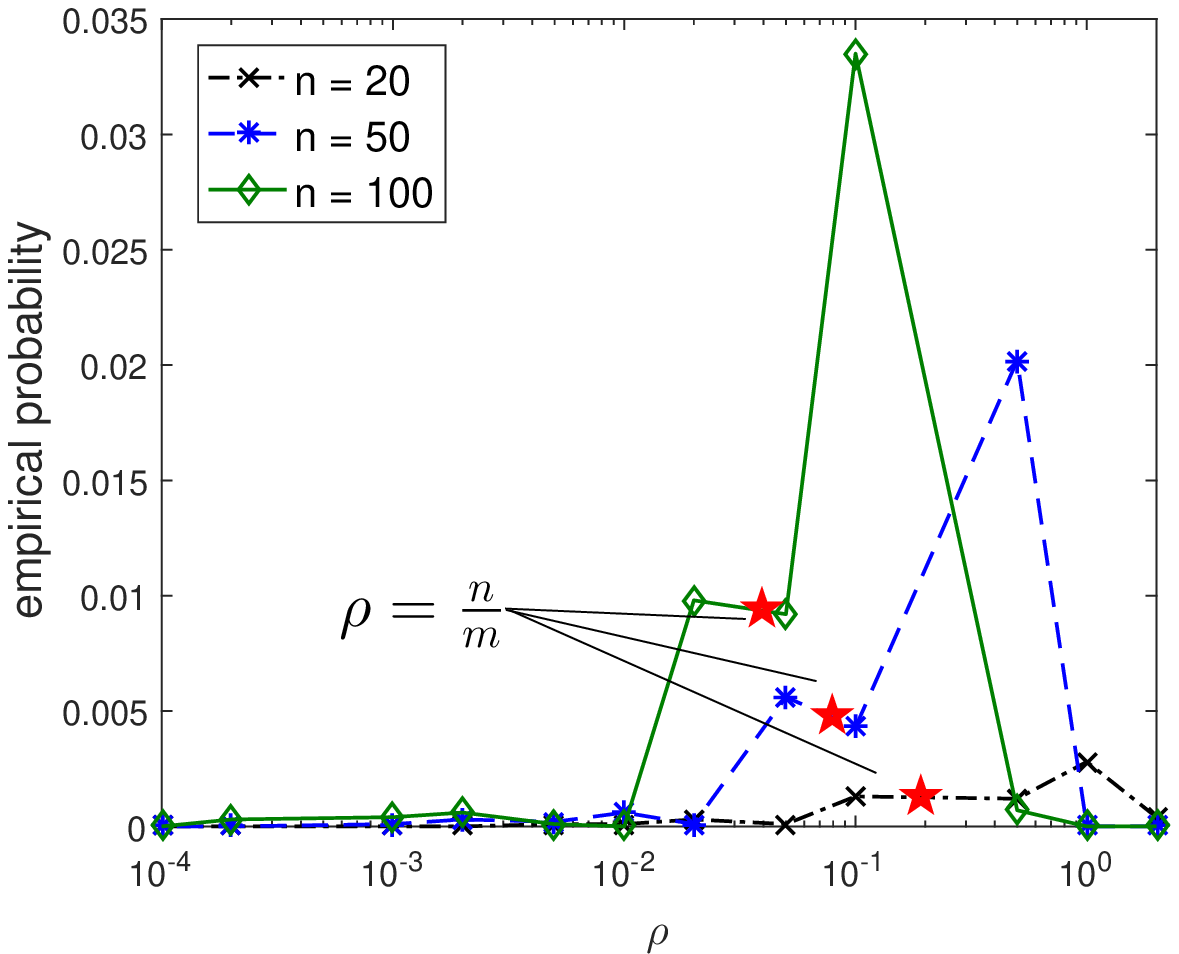}%
		\label{fig:CycHf}%
	}
	\subfigure[]{%
		\includegraphics[width=0.325\linewidth,height=0.26\linewidth]{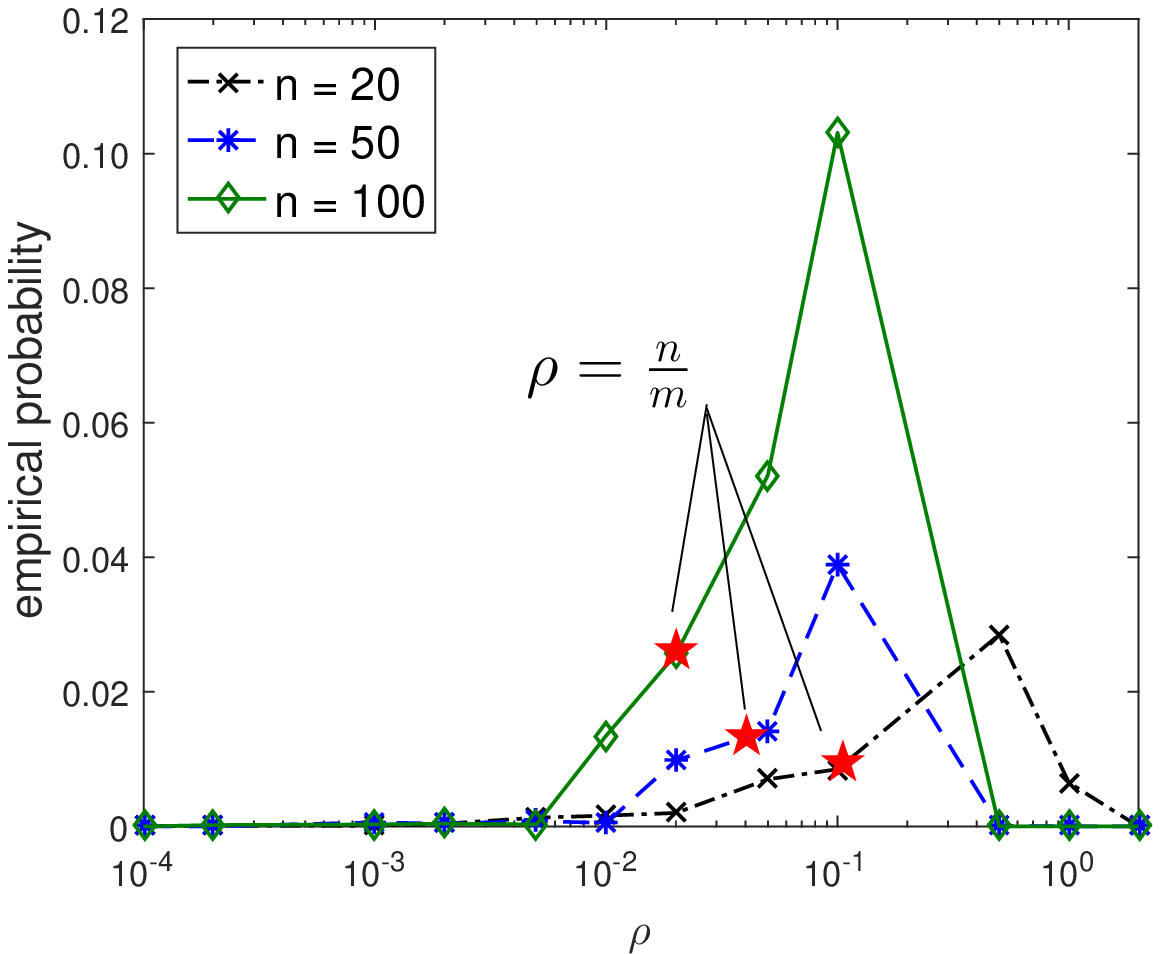}%
		\label{fig:CycCom}%
	}
	\caption{Empirical probability of cyclic case of BQ-CADMM in $10,000$ runs; $r_i\sim\mathcal{N}(0,100)+r_0$ with $r_0\sim\mathcal{N}(0,25)$, $L=30$ and $\Delta=1$. (a) star graph, (b) randomly generated graph with $m=\left\lceil\frac{(n+2)(n-1)}{4}\right\rceil$, (c) complete graph.}
	\label{fig:Cyc}
\end{figure*}

Since the convergent case can be viewed as the cyclic case with period $T=1$, we first investigate whether BQ-CADMM converges or cycles in various settings. Notice that the quantized consensus algorithm converges in all the above examples and indeed, BQ-CADMM and hence EBQ-CADMM tend to converge, particularly with small enough~$\rho$. To illustrate this, we simulate star graphs which have the smallest number of edges for a connected network, randomly generated connected graphs with intermediate numbers of edges, and complete graphs which have the highest number of edges. As Corollary~1 indicates that BQ-CADMM must converge when $\bar{r}$ is far from the bounded set $\mathcal{X}$, we set $L=30$ and generate $r_i\sim\mathcal{N}(0,100)+r_0$ where $r_0\sim\mathcal{N}(0,25)$ such that most simulated cases have $\bar{r}\in\mathcal{X}$. We then run BQ-CADMM with different $\rho$ for the same data. The simulation result is plotted in Fig.~\ref{fig:Cyc}, where each plotted value denotes the empirical probability of cyclic result in $10,000$ runs and both data and graph are randomly generated at each run.

Fig.~\ref{fig:Cyc} indicates that BQ-CADMM always converges for large enough $\rho$ with which the consensus result has a large error bound, as discussed in Section~\ref{sec:rhoCT}. When $\rho$ decreases, the number of cyclic cases first increases and then decreases to $0$. In particular, BQ-CADMM converges in all simulated examples with small enough $\rho$. Another interesting observation is that peak occurrence of the cyclic case differs for different network structures. The star network has highest cyclic cases around $\rho=0.5\frac{n}{m}$ and a larger network indicates less cyclic cases, while the intermediate and complete networks have cyclic cases centered around $\rho=5\frac{n}{m}$ with larger networks having more cyclic cases. 

In the same example, we also record the cyclic period when BQ-CADMM indeed cycles. We observe that the period of star networks is always $2$ for all $n$ and $\rho$, and the intermediate and complete networks have periods between $8$ and $15$. While we cannot draw firm conclusions on cyclic period from the simulation result, we do find that the period in all cyclic cases consists of two consecutive quantization levels for each node, which can help derive better error bounds as discussed in Section~\ref{sec:CP}.
\subsection{Convergence Time}
\label{sec:convergencetime}
To study how $\rho$ affects the convergence time of BQ-CADMM, we plot in Fig.~4 the average convergence time of the same example in the above section. We observe that BQ-CADMM converges immediately for large enough $\rho$ but again a large consensus error may exist. The convergence time is about $10$ when running BQ-CADMM with $\rho=\frac{n}{m}$ for all network sizes and network structures. For small $\rho$, the convergence time decreases polynomially as $\rho$ increases. The network structure also plays an important role here. For the star network, the convergence time is almost the same for all $n$ and  a larger $n$ results in a slightly longer convergence time. On the other hand, the intermediate and complete networks have longer convergence time with smaller $n$. Comparing the convergence time of the same network size, we also find that a denser network tends to converge faster for the same $\rho$.

\begin{figure*}[htbp]
	\centering
	\subfigure[]{%
		\includegraphics[width=0.325\linewidth,height=0.26\linewidth]{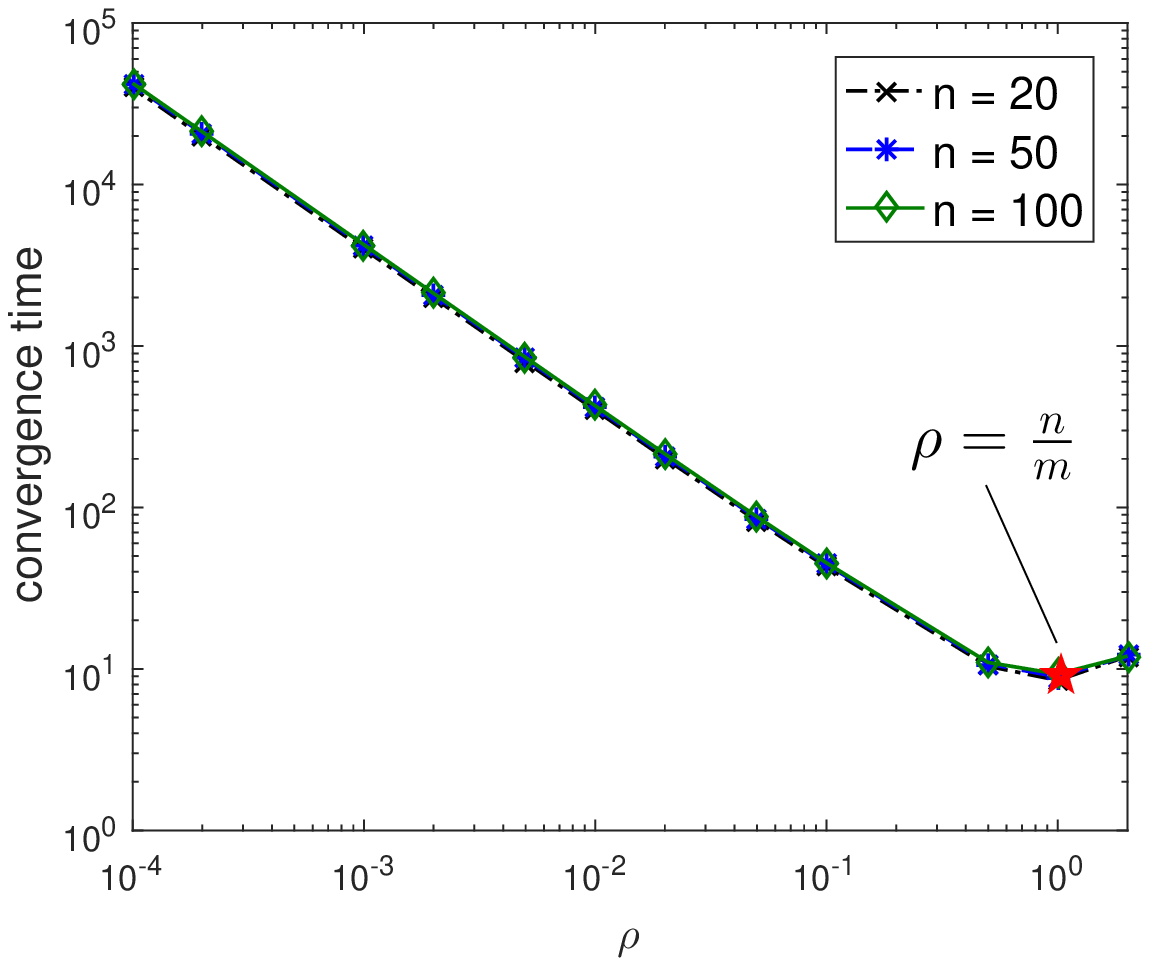}%
		\label{fig:CTStar}%
	}
	\subfigure[]{%
		\includegraphics[width=0.325\linewidth,height=0.26\linewidth]{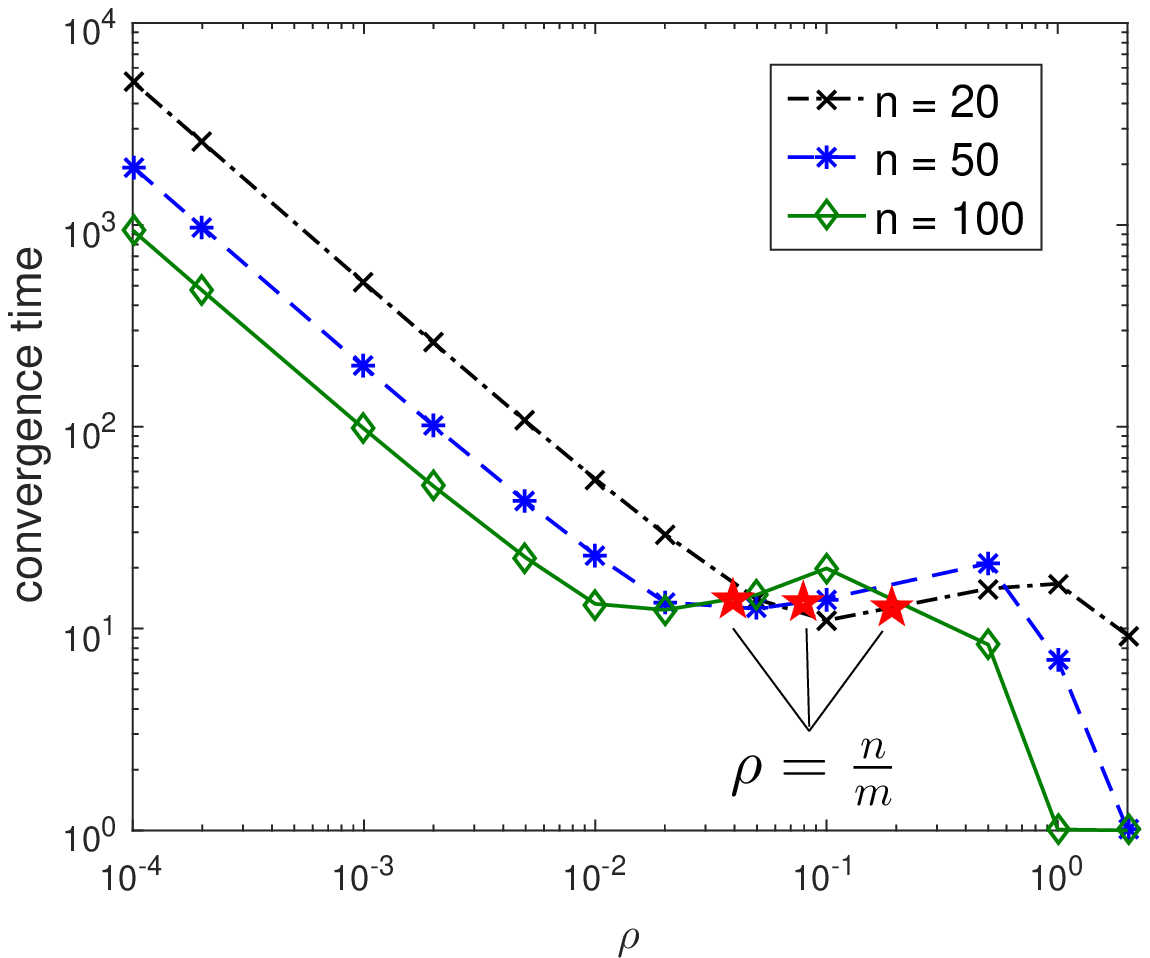}%
		\label{fig:CTHf}%
	}
	\subfigure[]{%
		\includegraphics[width=0.325\linewidth,height=0.26\linewidth]{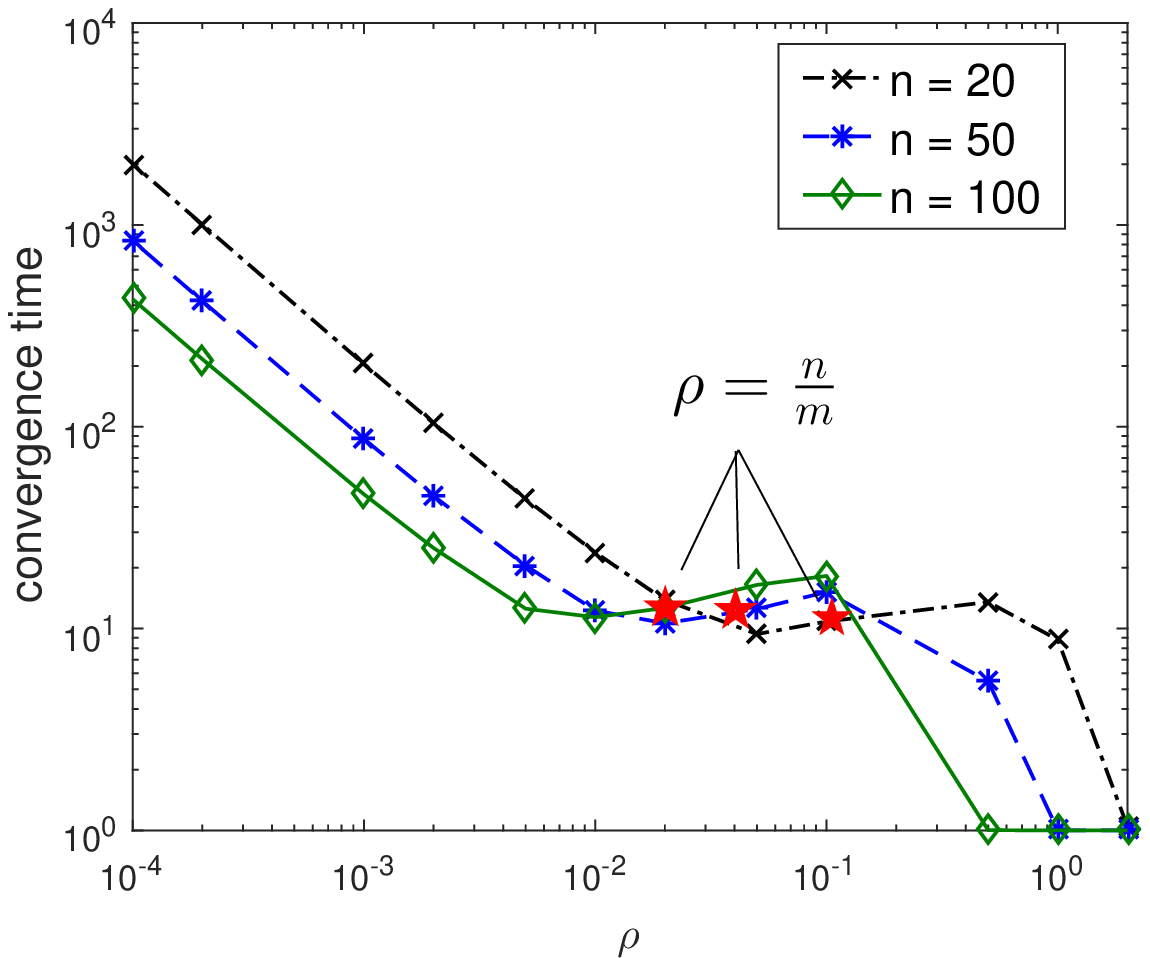}%
		\label{fig:CTCom}%
	}
	\caption{Convergence time of BQ-CADMM where each plotted value is the average of $10,000$ runs; $r_i\sim\mathcal{N}(0,100)+r_0$ with $r_0\sim\mathcal{N}(0,25)$, $L=30$ and $\Delta=1$. (a) star graph, (b) randomly generated graph with $m=\left\lceil\frac{(n+2)(n-1)}{4}\right\rceil$, (c) complete graph.}
	\label{fig:CT}
\end{figure*}

\subsection{Parameter Selection}
\label{sec:DSforPS}
To achieve high consensus accuracy (e.g., the consensus error is within one quantization resolution), the consensus error bound of the cyclic case in Theorem~\ref{thm:mainresults} implies that one may need to pick a very small $\rho$ which can make BQ-CADMM slow to reach convergence or cycling. Fortunately, simulations indicate that BQ-CADMM converges in most cases and $\rho$ can be larger to meet the consensus error requirement in convergent cases. As such, we can start BQ-CADMM with a large $\rho$ for a number of iterations. Even though the local variables may not converge, they become close to the optima. We then pick a small $\rho$ and run BQ-CADMM with the current local variable values. We continue this process until either a convergence is reached and $\rho$ is such that the consensus error bound (\ref{eqn:ProjCE}) meets the required consensus accuracy, or $\rho$ is small enough such that the cyclic error bound (\ref{eqn:cyclic_bd}) satisfies the accuracy requirement. We suggest starting with $\rho=\frac{n}{m}$ as it performs reasonably well in terms of both consensus accuracy and convergence time, as seen from  Fig.~\ref{fig:Cyc} and Fig.~\ref{fig:CT}. If only $n$ is known, we may simply start with $\rho=1$ as $n-1\leq m\leq\frac{n(n-1)}{2}$ for connected networks.

To further illustrate the strategy with decreasing step size, we next apply it to the same example in Section~\ref{sec:CP}. Starting with $\rho=\frac{n}{m}$, we run BQ-CADMM for $50$ iterations if $\rho >10^{-4}$ and then reduce it by a factor of $10$. We repeat this process until $\rho<10^{-4}$ with which we run BQ-CADMM long enough such that either convergence or cycling is reached. The average convergence time of this strategy is presented in Table~\ref{tab:comp}. As one can see, this strategy reduces dramatically the convergence time. In addition, most iterations occur when $\rho>10^{-4}$ and BQ-CADMM only takes a few more iterations to reach the convergence result in all simulated cases when $\rho\leq10^{-4}$. 

\begin{table}[ht]
\centering
\caption{Average convergence time of BQ-CADMM with and without using the decreasing strategy for parameter selection.}
{\renewcommand{\arraystretch}{1.25}
    \begin{tabular}{cccc}\hlinewd{.8pt}
\label{tab:comp}
        Structure           &~~Nodes~~   &Decreasing Parameter & Fixed Parameter \\
\hline
     \multirow{3}{*}{Star}                & $20$  & $251.2$               & $3.97\times 10^4$      \\ 
                    & $50$  & $253.9$               & $4.14\times 10^4$    \\ 
                      & $100$ & $257.1$               & $4.24\times 10^4$         \\ \hline
        \multirow{3}{*}{Intermediate}              & $20$  & $203.9$               & $5.13\times 10^3$          \\ 
            & $5$0  & $159.2$               & $1.94\times 10^3$         \\ 
                       & $100$ & $160.6$                     & $0.94\times10^3$        \\\hline
        \multirow{3}{*}{Complete}               & $20$  & $203.3$               & $2.00\times 10^3$ \\
                & $50$  & $158.7$               & $0.84\times 10^3$          \\ 
                       & $100$ & $160.4$               & $0.43\times 10^3$          \\\hline
        
    \end{tabular}}
\end{table}

\section{Conclusion and Discussion}
\label{sec:conclusion}
This paper proposes BQ-CADMM and EBQ-CADMM for distributed average consensus using finite-bit bounded quantization. It is shown that BQ-CADMM is equivalent to applying rounding quantizer to the CADMM update of a constrained least-squares problem. Based on this fact, we establish that BQ-CADMM either converges or cycles within finite iterations and further derive the consensus error bounds by utilizing convexity properties. An advantage of our algorithm is that the agents' data can be arbitrary; in contrast, existing quantized consensus algorithms that adopt finite-bit quantizers require the data to be bounded and a bound to be known. Noticing that the consensus error is large when the data average is much beyond the quantizer range, we proceed to develop EBQ-CADMM which is equivalent to repeated runs of BQ-CADMM. The resulting consensus error bound can be the same as that of using the unbounded rounding quantizer by picking a small enough algorithm parameter. Balancing fast convergence and small consensus error, an adaptive parameter selection scheme is devised that requires only the knowledge of the number of nodes. 

While simulations show that BQ-CADMM converges fast in most cases, the current work lacks theoretical analysis (e.g., upper bound) on cyclic period and convergence time. A possible approach is to find a Lyapunov function with which one can characterize these two quantities under suitable assumptions. In addition, asynchronous protocols, link failures, and noisy communication channels can be studied under the BQ-CADMM framework.

\appendix
\begin{IEEEproof}[Proof of Lemma~\ref{lem:alphabd}]
The proof is mainly based on the fact that $|\mathcal{Q}_b(\tilde{x})|\leq L$ for any $\tilde{x}\in\mathbb{R}$. We first prove the following:
\begin{align}
\mathcal{T_X}(x_{i}^{k+1})= 
\begin{cases}
-L,&{\text{if}}~\alpha_i^k>\left(1+4\rho|\mathcal{N}_i|\right)L+\left|r_i\right|,\nn\\
L, &{\text{if}}~\alpha_i^k<-\left(1+4\rho|\mathcal{N}_i|\right)L-\left|r_i\right|.\\
\end{cases}
\end{align}

Assume that $\mathcal{T_X}(x_{i}^{k+1}) >-L$ if $\alpha_{i}^k>\left(1+4\rho|\mathcal{N}_i|\right)L+\left|r_i\right|$. Similar to (\ref{eqn:ConstrainedDOpt}), the $x_i$-update of BQ-CADMM implies that $\mathcal{T_X}(x_{i}^{k+1})$ minimizes a constrained least-squares function:
\begin{align}
\mathcal{T_X}(x_{i}^{k+1})=&\hspace{2pt}\argmin_{\tilde{x}\in\mathcal{X}}\frac{1}{2}(\tilde{x}-r_i)^2+\rho|\mathcal{N}_i|\tilde{x}^2-\tilde{x}\Bigg(\rho|\mathcal{N}_i| \mathcal{Q}_b({x_i^k})\nn\\
&+\rho\sum_{j\in\mathcal{N}_i} \mathcal{Q}_b(x_{j}^k)-\alpha_i^k+r_i\Bigg)\nn\\
\triangleq&\hspace{2pt}\argmin_{\tilde{x}\in\mathcal{X}}G_i^k(\tilde{x}) + \alpha_i^k\tilde{x},\nn
\end{align}
where, for ease of presentation, we define
\begin{align}
G_i^k(\tilde{x})\triangleq&\hspace{2pt}\frac{1}{2}(\tilde{x}-r_i)^2+\rho|\mathcal{N}_i|\tilde{x}^2-\tilde{x}\Bigg(\rho|\mathcal{N}_i| \mathcal{Q}_b({x_i^k})\nn\\
&\hspace{2pt}+\rho\sum_{j\in\mathcal{N}_i} \mathcal{Q}_b(x_{j}^k)+r_i\Bigg).\nn
\end{align} 
Since $\left|\mathcal{Q}_b(\tilde{x})\right|\leq L$ for any $
\tilde{x}\in\mathbb{R}$, it is straightforward to verify that $G_i^k(\tilde{x})$ is Lipschitz continuous over $\mathcal{X}$: for any $\tilde{x},\tilde{y}\in\mathcal{X}$,
\begin{align}
\left|G_i^k(\tilde{x})-G_i^k(\tilde{y})\right|<\left(\left(1+4\rho|\mathcal{N}_i|\right)L+\left|r_i\right|\right)\left|\tilde{x}-\tilde{y}\right|.\nn
\end{align}
Now letting $\tilde{x}=-L\in\mathcal{X}$, we have
\begin{align}
&\hspace{2pt}G_i^k(\tilde{x}) + \alpha_i^k\tilde{x}-\left(G_i^k\left(\mathcal{T_X}(x_i^{k+1})\right)+ \alpha_i^k\mathcal{T_X}(x_i^{k+1})\right)\nn\\
=&\hspace{2pt}G_i^k(-L) -G_i^k\left(\mathcal{T_X}(x_i^{k+1})\right)- \alpha_i^kL-\alpha_i^k\mathcal{T_X}(x_i^{k+1})\nn\\
<&\hspace{2pt}\left(\left(1+4\rho|\mathcal{N}_i|\right)\hspace{-1pt}L+\left|r_i\right|-\alpha_i^k\right)\hspace{-1pt}\left(\mathcal{T_X}(x_i^{k+1})+L\right)\hspace{-1pt},\nn\\
<&\hspace{2pt}0,\nn
\end{align}
where the last inequality is because $\alpha_i^k>\left(1+4\rho|\mathcal{N}_i|\right)L+\left|r_i\right|$ and $\mathcal{T_X}(x_i^{k+1})>-L$. This contradicts the fact that $\mathcal{T_X}(x_i^{k+1})$ minimizes $G_i^k(\tilde{x}) + \alpha_i^k\tilde{x}$ over  $\mathcal{X}$. That $\mathcal{T_X}(x_{i}^{k+1})=L$ for $\alpha_i^k<-\left(1+4\rho|\mathcal{N}_i|\right)L-\left|r_i\right|$ can be shown analogously.

If $\mathcal{T_X}(x_{i}^{k+1})=-L$, the bounded quantization scheme implies $\mathcal{Q}_b(x_i^{k+1})=-L\leq\mathcal{Q}_b(x_j^{k+1})$ for any $j\in\mathcal{N}_i$. Therefore, if $\alpha_i^k>(1+4\rho|\mathcal{N}_i|)L+|r_i|$, we get
\begin{align}
\alpha_{i}^{k+1}=\alpha_{i}^k+\rho|\mathcal{N}_i|\mathcal{Q}_b(x_i^{k+1})-\rho\sum_{j\in\mathcal{N}_i}\mathcal{Q}_b(x_j^{k+1})\leq\alpha_{i}^k.\nn
\end{align}
Similarly, $\alpha_{i}^{k+1}\geq\alpha_{i}^k$ if $\alpha_{i}^k<-(1+4\rho|\mathcal{N}_i|)L-|r_i|$. 

Next consider $|\alpha_{i}^k|\leq (1+4\rho|\mathcal{N}_i|)L+|r_i|$. In this case, we can simply use the triangle inequality to conclude that
\begin{align}
|\alpha_{i}^{k+1}|&\leq|\alpha_{i}^k|+\rho|\mathcal{N}_i|\left|\mathcal{Q}_b(x_{i}^{k+1})\right|+\rho\sum_{j\in\mathcal{N}_i}\left|\mathcal{Q}_b(x_{i}^{k+1})\right|\nn\\
&\leq(1+6\rho|\mathcal{N}_i|)L+|r_i|.\nn
\end{align}

Since $\alpha_i^0=0$, we finally have $|\alpha_{i}^k|\leq(1+6\rho|\mathcal{N}_i|)L+|r_i|$ for all $k$.
\end{IEEEproof}
\bibliographystyle{IEEEtran}
\bibliography{SZhuBib}

\end{document}